\documentclass{amsart}
\usepackage{amsmath,amssymb,amsxtra,amsthm,amscd}
\usepackage{lmodern}
\usepackage[all,cmtip]{xy}
\usepackage{rotating}
\usepackage{microtype}
\usepackage{todonotes}

\usepackage{adjustbox}
\usepackage{color}

\newif\ifShowLabels
\ShowLabelstrue
\newcommand{\TeXref}[1]{
\marginpar{\scriptsize \texttt{#1}}}

\ShowLabelsfalse

\DeclareMathOperator{\EGamma}{\boldsymbol{E}\boldsymbol{\Gamma}}

\DeclareMathOperator{\fil}{fil}

\DeclareMathOperator{\Free}{\mathbf{Free}}

\DeclareMathOperator{\Fun}{Fun}

\DeclareMathOperator{\Hom}{Hom}

\DeclareMathOperator{\id}{id}

\DeclareMathOperator{\K}{\mathit{K}}
         \newcommand{\Knc}{\K^{-\infty}}

         \newcommand{\KncX}{\K^{-\infty}_X}

\DeclareMathOperator{\point}{pt}

\DeclareMathOperator{\Sh}{Sh}

\DeclareMathOperator{\W}{\mathcal{W}}

\DeclareMathOperator*{\one}{1}
\newcommand{\onehatplace}[1]
{ \one^{\substack{#1 \\ \frown}} }

\DeclareMathOperator*{\bones}{\times}
\newcommand{\undertimes}[1]
{ \bones_{#1} }

\DeclareMathOperator*{\bowl}{\cup}
\newcommand{\undercup}[1]
{ \bowl_{#1} }

\DeclareMathOperator*{\arch}{\cap}
\newcommand{\undercap}[1]
{ \arch_{#1} }

\newcommand{\pull}
{\!\!\! -\!\!\! -\!\!\! -\!\!\!}

\DeclareMathOperator*{\holimprep}{holim}                       
\newcommand{\holim}[1]%
{\displaystyle\holimprep_{\substack{\leftarrow \pull - \\ #1}} \, }

\DeclareMathOperator*{\hocolimprep}{hocolim}                   
\newcommand{\hocolim}[1]%
{\displaystyle\hocolimprep_{\substack{- \pull \rightarrow \\ #1}} \, }

\DeclareMathOperator*{\plainlim}{lim}                           
\newcommand{\contralim}[1]%
{\displaystyle\plainlim_{\substack{\leftarrow \pull - \\ #1}} \, }

\DeclareMathOperator*{\plaincolim}{colim}                       
\newcommand{\colim}[1]%
{\displaystyle\plaincolim_{\substack{- \pull \rightarrow \\ #1}} \, }

\DeclareMathOperator*{\laxlimplain}{laxlim}                     
\newcommand{\laxlim}[1]%
{\displaystyle\laxlimplain_{\substack{\leftarrow \pull - \\ #1}} \, }

\providecommand{\bysame}{\makebox[3em]{\hrulefill}\thinspace}

\theoremstyle{plain}
\newtheorem{Thm}{Theorem}[subsection]

\newtheorem{Cor}[Thm]{Corollary}

\newtheorem{Lem}[Thm]{Lemma}
\newtheorem{Prop}[Thm]{Proposition}

\theoremstyle{definition}
\newtheorem{Def}[Thm]{Definition}

\newtheorem{Ex}[Thm]{Example}
\newtheorem{Exs}[Thm]{Examples}

\newtheorem{Rem}[Thm]{Remark}

\newtheorem{Que}[Thm]{Question}
\newtheorem{Conj}[Thm]{Conjecture}

\theoremstyle{remark}
\newtheorem{Not}[Thm]{Notation}

\newtheoremstyle{freestylethm}{6pt}{6pt}{\itshape}{}%
                {\bfseries}{}{.5em}{\thmnote{#3}}
\theoremstyle{freestylethm}

\newcommand{\SecRef}[2]{\section{#1}\label{S:#2}%
\ifShowLabels \TeXref{{S:#2}} \fi}
\newcommand{\SSecRef}[2]{\subsection{#1}\label{SS:#2}%
\ifShowLabels \TeXref{{SS:#2}} \fi}

\newcommand{\refS}[1]{\textup{\ref{S:#1}}}
\newcommand{\refSS}[1]{\textup{\ref{SS:#1}}}

\newcommand{\refT}[1]{\textup{\ref{T:#1}}}
\newcommand{\refL}[1]{\textup{\ref{L:#1}}}
\newcommand{\refD}[1]{\textup{\ref{D:#1}}}
\newcommand{\refC}[1]{\textup{\ref{C:#1}}}

\newcommand{\refP}[1]{\textup{\ref{P:#1}}}

\newcommand{\refN}[1]{\textup{\ref{N:#1}}}
\newcommand{\refQ}[1]{\textup{\ref{Q:#1}}}

\newenvironment{ThmRef}[1]%
{ \begin{Thm} \label{T:#1}
\ifShowLabels \TeXref{T:#1} \fi }%
{ \end{Thm} }
\newenvironment{DefRef}[1]%
{ \begin{Def} \label{D:#1}
\ifShowLabels \TeXref{D:#1} \fi }%
{ \end{Def} }
\newenvironment{LemRef}[1]%
{ \begin{Lem} \label{L:#1}
\ifShowLabels \TeXref{L:#1} \fi }%
{ \end{Lem} }
\newenvironment{CorRef}[1]%
{ \begin{Cor} \label{C:#1}
\ifShowLabels \TeXref{C:#1} \fi }%
{ \end{Cor} }
\newenvironment{RemRef}[1]%
{ \begin{Rem} \label{R:#1}
\ifShowLabels \TeXref{R:#1} \fi }%
{ \end{Rem} }
\newenvironment{PropRef}[1]%
{ \begin{Prop} \label{P:#1}
\ifShowLabels \TeXref{P:#1} \fi }%
{ \end{Prop} }
{ \begin{Ex} \label{E:#1}
\ifShowLabels \TeXref{E:#1} \fi  }%
{ \end{Ex} }
\newenvironment{NotRef}[1]%
{ \begin{Not} \label{N:#1}
\ifShowLabels \TeXref{N:#1} \fi }%
{ \end{Not} }
\newenvironment{QueRef}[1]%
{ \begin{Que} \label{Q:#1}
\ifShowLabels \TeXref{Q:#1} \fi }%
{ \end{Que} }
\newenvironment{ConjRef}[1]%
{ \begin{Conj} \label{Conj:#1}
\ifShowLabels \TeXref{Conj:#1} \fi }%
{ \end{Conj} }

\newenvironment{ThmRefName}[2]%
{ \begin{Thm} [#2]\label{T:#1}
\ifShowLabels \TeXref{T:#1} \fi }%
{ \end{Thm} }
\newenvironment{DefRefName}[2]%
{ \begin{Def} [#2]\label{D:#1}
\ifShowLabels \TeXref{D:#1} \fi }%
{ \end{Def} }
{ \begin{Lem} [#2]\label{L:#1}
\ifShowLabels \TeXref{L:#1} \fi }%
{ \end{Lem} }
{ \begin{Cor} [#2]\label{C:#1}
\ifShowLabels \TeXref{C:#1} \fi }%
{ \end{Cor} }
{ \begin{Rem} [#2]\label{R:#1}
\ifShowLabels \TeXref{R:#1} \fi }%
{ \end{Rem} }
{ \begin{Prop} [#2]\label{P:#1}
\ifShowLabels \TeXref{P:#1} \fi }%
{ \end{Prop} }
{ \begin{Ex} [#2]\label{E:#1}
\ifShowLabels \TeXref{E:#1} \fi }%
{ \end{Ex} }

\let\oldtocsection=\tocsection
\let\oldtocsubsection=\tocsubsection
\renewcommand{\tocsection}[2]{\hspace{0em}\oldtocsection{#1}{#2}}
\renewcommand{\tocsubsection}[2]{\hspace{2em}\oldtocsubsection{#1}{#2}}

\begin{document}

\title[\textit{K}-theory with fibred control]{\textit{K}-theory of geometric modules with fibred control}
\author[Gunnar Carlsson]{Gunnar Carlsson}
\address{Department of Mathematics\\ Stanford University\\ Stanford\\ CA 94305}
\email{gunnar@math.stanford.edu}
\author[Boris Goldfarb]{Boris Goldfarb}
\address{Department of Mathematics and Statistics\\ SUNY\\ Albany\\ NY 12222}
\email{goldfarb@math.albany.edu}
\date{\today}

\begin{abstract}
Controlled algebra plays a central role in many recent advances in geometric topology.
This paper studies the iteration construction that was present from the very origins of the theory but started being exploited only recently.
We develop the general framework for fibred control, prove localization theorems required for fibrewise excision, and then prove several versions of fibrewise excision theorems. However, we also demonstrate how the standard tools break down in the presence of new equivariant phenomena which require advanced localization methods.
\end{abstract}

\maketitle

\tableofcontents

\SecRef{Introduction}{INTRO}

Controlled algebra  of geometric modules has been a very useful device in geometric topology since its appearance in the work of Frank Quinn. The most effective work on Novikov, Borel, Farrell-Jones conjectures and other prominent problems uses geometric control conditions well-suited for generating $K$-theory spectra out of the associated additive controlled categories of free modules.

The bounded controlled algebra was developed by E.K. Pedersen and C. Weibel \cite{ePcW:85,ePcW:89}.  It has been used extensively for delooping $K$-theory of rings and other algebraic objects and in the study of the Novikov conjecture \cite{gC:95,gCbG:04,gCeP:93} about the assembly map for algebraic $K$-theory of certain group rings. 
   
We start by explaining this theory as an algebraic theory of free modules, possibly infinitely generated, parametrized by a metric space.  The input is a proper metric space $X$ (i.e. a metric space in which every closed bounded subset is compact) and a ring $R$.  Given this data, consider triples $(M, B, \varphi)$, where $M$ is a free left $R$-module with basis $B$, and where $\varphi \colon B \rightarrow X$ is a reference function with the property that the inverse image of any bounded subset is finite.  A morphism from $(M,B,\varphi)$ to $(M^{\prime}, B^{\prime} , \varphi^{\prime})$  is an $R$-module homomorphism $f$ from $M$ to $M^{\prime}$  which has the property that there exists a bound $b \ge 0$ so that for any basis element $\beta \in B$, $f(\beta)$ is in the span of basis elements $\beta' \in B^{\prime}$ for which $\varphi' (\beta ^{\prime}) \leq b$.  This is the simplest version of a {\em control condition} that can be imposed on homomorphisms to construct various categories of modules, and it leads to the category of \textit{geometric modules} $\mathcal{C} (X, R)$.  This is an additive category to which one can apply the usual algebraic $K$-theory construction.   

A formulation of this theory that was observed already in Pedersen-Weibel \cite{ePcW:85} allows one to use objects from an arbitrary additive category $\mathcal{A}$ as ``coefficients'' in place of the finitely generated free $R$-modules.  We will spell out the details of this formulation further in the paper.  So one obtains a new additive category $\mathcal{C} (X, \mathcal{A})$. This observation leads to the possibility of iterating the geometric control construction: if there are two proper metric spaces $X$ and $Y$, then we have an additive category $\mathcal{C} (X, \mathcal{C} (Y, R))$.  This is precisely a case of \textit{fibred control} we want to study in this paper. To give an idea for what kind of control is implied by this construction, let's parse the implications for objects viewed as $R$-modules.  In this setting an object is a free module $M$ viewed as parametrized over the product $X \times Y$, so we have a pair $(M,B)$ with a reference function $\varphi \colon B \rightarrow X \times Y$.  The control condition on homomorphisms $f$ from $(M_1, B_1, \varphi _1)$ to $(M_2, B_2, \varphi _2)$ amounts to the existence of a number $b_X \ge 0$ and a function $c \colon X \to [0, \infty)$ such that for $\beta \in B_1$, we have that $f(\beta)$ is a linear combination of basis elements $\beta ^{\prime} \in B_2$ so that $d(\pi _X(\varphi _1(\beta)), \pi _X(\varphi _2 (\beta ^{\prime}))) \leq b_X$ and 
 $d(\pi_Y (\varphi _1(\beta)), \pi_Y (\varphi _2 (\beta ^{\prime}))) \leq c (\pi _X(\varphi _1(\beta)))$.  
Contrast this with the control condition in $\mathcal{C} (X \times Y, R)$ which is exactly as above except $c$ can be chosen to be a constant function.
This should suggest that the condition of fibred control, with $X$ regarded as a base and $Y$ regarded as a fibre, relaxes the usual control condition over $X \times Y$ by letting $c$ be a quantity varying with $x \in X$. 

It is shown in \cite{gC:95} that $\mathcal{C} (X, R)$ enjoys a certain excision property, which permits among other things a comparison with Borel-Moore homology spectra with coefficients in the $K$-theory spectrum $K(A)$. 
This property becomes crucial in the equivariant homotopy theoretic approach to the Novikov conjecture.
In this paper we prove a variety of excision results in categories with fibred control and some related categories of fixed-point objects with respect to naturally occurring actions on the metric spaces.
These theorems are required for further work on new cases of the Novikov conjecture and other conjectures about assembly maps in $K$-theory.

\medskip

The generalization of controlled algebra in this paper 
is useful for a variety of applications.  

It is most immediately the natural controlled theory to consider for bundles on non-compact manifolds where the control in the fibre direction can vary with the fibre.
This holds more generally for stacks or the coarse analogues of bundles that have already appeared in the literature \cite{kW:10}. 
There are also important constructions that are popular in coarse geometry as extensions in so-called Fibering Permanence theorems.  They can be thought of as total spaces of fibration-like maps where instead of local trivializations one postulates a specific coarse property $P$ that preimages of metric balls need to satisfy. The permanence theorems then derive the same property $P$ for the total space. 
Now suppose some geometric property $P$ allows a bounded excision strategy as developed in \cite{gC:95,gCbG:04} and reviewed in section \refSS{BETRGCV}.  From our fibred excision theorems, one learns how to compute bounded $K$-theory of the total space.  This becomes most interesting when the permanence actually fails, and the total space is not known to satisfy the same property $P$.  The reason is that even though the space does not satisfy property $P$, we have expanded the class of spaces with computatble bounded $K$-theory.  Examples of such situations and related discussion can be found in \cite{sBbG:18,bGjG:18}.

Another manifestation of fibred control is found in the active current research area seeking sufficient geometric conditions for verifying the coarse Baum-Connes conjecture.  The original condition in \cite{gY:98,gY:00}  was existence of a coarse embedding of the group into a Hilbert space. The embedding can be viewed as a way to gain geometric control for performing stable excision over the group or the space with a cocompact action by the group.   This condition has been relaxed by several authors to \textit{fibred coarse embeddings} \cite{sA:16,xCqWxW:13,xCqWgY:13,mF:14,mMhS:18}.  The theorems in this paper allow to prove the  coarse Borel conjecture in $K$-theory for the same class of spaces.

\medskip

Finally, let us address the methods used and the interesting phenomena that come up.  The go-to technique for proving localization and excision theorems in controlled algebra is the use of Karoubi filtrations in additive categories.  This technique is indeed what is used also in our proofs of general non-equivariant fibrewise excision results (section \refSS{FETDSLK}) and in some specific equivariant situations (section \refSS{FLEKth}).  However, we want point out one curious situation that will interest the experts. The most useful novelty of the fibred control is that it allows to introduce constraints on features, including actions, that vary across fibers.  When this is done to actions, it becomes unnatural to insist that the action is by maps that are necessarily isometries on the nose or, if the action is by more general coarse equivalences, that the same numerical constraint is satisfied across all fibers.  It turns out that all of this leads to some desired excision statements about fixed point object categories for multiple actions that might be true but are inaccessible to the Karoubi filtration technique.  We give an example of such situation in the last section of the paper where we are able to pinpoint the deficiency that makes the Karoubi filtrations unavailable. 

There are two viable options for dealing with this problem.  Both of them are developed by the authors in separate papers.  One leads to the development of fibred $G$-theory as part of controlled $G$-theory \cite{gCbG:00,gCbG:15} which is known to have better localization properties.  There is also a comparison Cartan map from $K$-theory to $G$-theory and a set of fairly general conditions under which the Cartan map is an equivalence.
The other option is a different larger context that we call approximate $K$-theory.  This second theory has additional excision properties which turn out to be essential for the $L$-theory version of this material.

\SecRef{Pedersen-Weibel categories, Karoubi filtrations}{KFPWC}

\SSecRef{Pedersen-Weibel theory}{PWTH}
Bounded $K$-theory introduced in Pedersen \cite{eP:84} and Pedersen--Weibel \cite{ePcW:85} associates a nonconnective spectrum $\Knc (M,R)$ to a proper metric space $M$ (a metric space where closed bounded subsets are compact) and an associative ring $R$ with unity.  We are going to start with a careful, at times revisionist, review of the well-known features of this theory that will need generalization.   

\begin{DefRef}{PWcats}
$\mathcal{C} (M,R)$ is the additive category of \textit{geometric $R$-modules} whose 
objects are functions $F \colon M \to \Free_{fg} (R)$ which are locally finite assignments of free finitely generated $R$-modules $F_m$ to points $m$ of $M$.
The local finiteness condition requires precisely that for any bounded subset $S \subset M$ the restriction of $F$ to $S$ has finitely many nonzero values.
Let $d$ be the distance function in $M$.  The morphisms in $\mathcal{C} (M,R)$ are the $R$-linear homomorphisms
\[
\phi \colon \bigoplus_{m \in M} F_m \longrightarrow \bigoplus_{n \in M} G_n
\]
with the property
that the components $F_m \to G_n$ are zero for $d(m,n) > b$
for some fixed real number $b = b (\phi) \ge 0$.
The associated $K$-theory spectrum is denoted by $K (M,R)$, or $K (M)$ when the choice of ring $R$ is implicit,
and is called the \textit{bounded K-theory} of $M$.
\end{DefRef}

\begin{NotRef}{ttyh}
For a subset $S \subset M$ and a real number $r \ge 0$,
$S[r]$ will stand for the metric $r$-enlargement $\{ m \in M \mid d (m,S) \le r \}$.
In this notation, the metric ball of radius $r$ centered at $x$ is $\{ m \} [r]$ or simply $m [r]$.
\end{NotRef}

Now for every object $F$ there is a free $R$-module $F(S) = \bigoplus_{m \in S} F_m$.  The condition that $\phi$ is controlled as above is equivalent to existence of a number $b \ge 0$ so that $\phi F(S) \subset F(S[b])$ for all choices of $S$.

\begin{PropRef}{KLOWAS}
	The description of $\mathcal{C} (M,R)$ in the Introduction defines a category additively equivalent to the one in Definition \refD{PWcats}, establishing a dictionary between terminology in various papers in the literature.
\end{PropRef}

\begin{proof}
	Given a triple $(M, B, \varphi)$ described in the Introduction, define $M_x$ as the $R$-submodule freely generated by $\varphi^{-1} (x)$.  It is clear that this gives an additive functor in one direction.  The inverse functor is constructed by selecting a finite basis in each $F_x$ and defining $B$ to be the union of these bases.  The map $\varphi$ sends a basis element $b$ to $x$ if $b \in F_x$.
\end{proof}

Inclusions of metric spaces induce additive functors between the corresponding bounded $K$-theory spectra.
The main result of Pedersen--Weibel \cite{ePcW:85} is a delooping theorem which can be stated as follows.

\begin{ThmRefName}{Deloop}{Nonconnective delooping of bounded \textit{K}-theory}
Given a proper metric space $M$ and the standard Euclidean metric on the real line $\mathbb{R}$, the natural inclusion $M \to M \times \mathbb{R}$ induces isomorphisms $K_n (M) \simeq K_{n-1} (M \times \mathbb{R})$
for all integers $n > 1$.
If one defines the spectrum
\[
\Knc (M,R) \ = \ \hocolim{k} \Omega^k K (M \times \mathbb{R}^k),
\]
then the stable homotopy groups of  $\Knc (R) = \Knc (\point,R)$ coincide with the algebraic $K$-groups of $R$ in positive dimensions and with the Bass negative $K$-theory of $R$ in negative dimensions.
\end{ThmRefName}

\SSecRef{Bounded excision theorem}{BETRGCV}

Suppose $U$ is a subset of $M$.
Let $\mathcal{C} (M,R)_{<U}$ denote the full subcategory of $\mathcal{C} (M,R)$ on the objects $F$ with $F_m = 0$ for all points $m \in M$ with $d (m, U) > D$ for some fixed number $D > 0$ specific to $F$.
This is an additive subcategory of $\mathcal{C} (M,R)$ with the associated $K$-theory spectrum $\Knc (M,R)_{<U}$.
Similarly, if $U$ and $V$ are a pair of subsets of $M$, then there is the full additive subcategory
$\mathcal{C} (M,R)_{<U,V}$ of $F$ with $F_m = 0$ for all $m$ with $d (m, U) \le D_1$ and $d (m, V) \le D_2$ for some numbers $D_1, D_2 > 0$.
It is easy to see that $\mathcal{C} (M,R)_{<U}$ is in fact equivalent to $\mathcal{C} (U,R)$.

\medskip

The following theorem is the basic computational device in bounded $K$-theory.

\begin{ThmRefName}{BddExc}{Bounded excision \cite{gC:95}}
Given a proper metric space $M$ and a pair of subsets $U$, $V$ of $M$, there is a homotopy pushout diagram
\[
\xymatrix{
 \Knc (M)_{<U, V} \ar[r] \ar[d]
&\Knc (U) \ar[d] \\
 \Knc (V) \ar[r]
&\Knc (M) }
\]
\end{ThmRefName}

It is possible to restate the Bounded Excision Theorem in a more intrinsic form, after restricting to a special class of coverings.

\begin{DefRef}{CATHJK}
A pair of subsets $S$, $T$ of a proper metric space $M$ is called \textit{coarsely antithetic} if $S$ and $T$ are proper metric subspaces with the subspace metric and for each number $K > 0$ there is a number $K' > 0$ so that
\[
S[K] \cap T[K] \subset (S \cap T)[K'].
\]
\end{DefRef}

Examples of coarsely antithetic pairs
include any two non-vacuously intersecting closed subsets of a simplicial tree as well as
complementary closed half-spaces in a Euclidean space.

\begin{CorRef}{BDDEXCII}
If $U$ and $V$ is a coarsely antithetic pair of subsets of $M$ which form a cover of $M$, then the commutative square
\[
\xymatrix{
 \Knc (U \cap V) \ar[r] \ar[d]
&\Knc (U) \ar[d] \\
 \Knc (V) \ar[r]
&\Knc (M)
}
\]
is a homotopy pushout.
\end{CorRef}

We want to outline the proof of Theorem \refT{BddExc} in specific terms that will be used later.

The notion of Karoubi filtrations in additive categories is central to the proof of this theorem.
The details can be found in Cardenas--Pedersen \cite{mCeP:97}.

\begin{DefRef}{Karoubi}
An additive category $\mathcal{C}$ is \textit{Karoubi filtered} by a full subcategory $\mathcal{A}$
if every object $C$ has a family of decompositions $\{ C = E_{\alpha} \oplus D_{\alpha} \}$ with $E_{\alpha} \in \mathcal{A}$ and $D_{\alpha} \in \mathcal{C}$, called a \textit{Karoubi filtration} of $C$, satisfying the following properties.
\begin{itemize}
\item For each $C$, there is a partial order on Karoubi decompositions such that $E_{\alpha} \oplus D_{\alpha} \le E_{\beta} \oplus D_{\beta}$ whenever $D_{\beta} \subset D_{\alpha}$ and $E_{\alpha} \subset E_{\beta}$.
\item Every map $A \to C$ factors as $A \to  E_{\alpha} \to E_{\alpha} \oplus D_{\alpha} = C$ for some $\alpha$.
\item Every map $C \to A$ factors as $C = E_{\alpha} \oplus D_{\alpha} \to  E_{\alpha} \to A$ for some $\alpha$.
\item For each pair of objects $C$ and $C'$ with the corresponding filtrations $\{ E_{\alpha} \oplus D_{\alpha} \}$ and $\{ E'_{\alpha} \oplus D'_{\alpha} \}$, the filtration of $C \oplus C'$ is the family $\{ C \oplus C' = ( E_{\alpha} \oplus E'_{\alpha} ) \oplus ( D_{\alpha} \oplus D'_{\alpha} ) \}$.
\end{itemize}
A morphism $f \colon C \to D$ is $\mathcal{A}$-\textit{zero} if $f$ factors through an object of $\mathcal{A}$.
Define the \textit{Karoubi quotient} $\mathcal{C}/\mathcal{A}$ to be the additive category with the same objects as $\mathcal{C}$ and morphism sets $\Hom_{\mathcal{C}/\mathcal{A}} (C,D) = \Hom (C,D)/ \{ \mathcal{A}\mathrm{-zero \ morphisms} \}$.
\end{DefRef}

The following is the main theorem of Cardenas--Pedersen \cite{mCeP:97}.

\begin{ThmRefName}{CPF}{Fibration theorem}
Suppose $\mathcal{C}$ is an $\mathcal{A}$-filtered category, then there is a homotopy fibration
\[
K ( \mathcal{A}^{\wedge K} ) \longrightarrow
K ( \mathcal{C} ) \longrightarrow K ( \mathcal{C}/\mathcal{A} ).
\]
Here $\mathcal{A}^{\wedge K}$ is a certain subcategory of the idempotent completion of $\mathcal{A}$ with the same positive $K$-theory as $\mathcal{A}$.
\end{ThmRefName}

By an observation in \cite{gC:95}, one immediately gets the following consequence.

\begin{CorRef}{CPFrev}
Suppose $\mathcal{C}$ is an $\mathcal{A}$-filtered category, then there is a homotopy fibration
\[
\Knc ( \mathcal{A} ) \longrightarrow
\Knc ( \mathcal{C} ) \longrightarrow \Knc ( \mathcal{C}/\mathcal{A} ).
\]
\end{CorRef}

Now Theorem \refT{BddExc} follows from this Corollary by the following device.

The first crucial observation is that $\mathcal{C} = \mathcal{C}(M)$ is $\mathcal{A} = \mathcal{C} (M)_{< U}$-filtered.
For a future reference, let us spell out what is involved.
The additive structure in $\mathcal{C}$ is given by $(F \oplus G)_m = F_m \oplus G_m$.  So, in particular, $(F \oplus G)(S) = F(S) \oplus G(S)$ for all subsets $S$.
Given an object $F$ of $\mathcal{C}$, the subobjects $F (U[k])$, $k \ge 0$, are free direct summands of $F(M)$ which can be given the structure of a geometric module over $M$ in the obvious way.  Now the decompositions $F = F (U[k]) \oplus F (M \setminus U[k])$ is the family we need.  Suppose, for illustration, we have $f \colon A \to F$ bounded by $b \ge 0$.  Then 
$A(M) = A(U[r])$ for some $r \ge 0$ so $f (A) \subset F(U[r + b])$.  So indeed $f$ factors through this direct summand.

Since $\mathcal{C}(M)$ is $\mathcal{C} (M)_{< V}$-filtered, there is the additive Karoubi quotient which we denote by $\mathcal{C} (M,V)$, with the nonconnective $K$-theory $\Knc (M, V)$.  For simplicity, let us assume that $U$ and $V$ form a coarsely antithetic pair, then $\mathcal{C} (M)_{< U}$ is similarly $\mathcal{C} (M)_{< U,V}$-filtered, with the Karoubi quotient $\mathcal{C} (U, U \cap V)$.
Corollary \refC{CPFrev} gives two homotopy fibrations that form a commotative diagram

\[
\begin{CD}
\Knc (M)_{< U,V} @>>> \Knc (M)_{< U} @>>> \Knc (U, U \cap V) \\
@VVV @VVV @VV{\simeq}V \\
\Knc (M)_{< V} @>>> \Knc (M) @>>> \Knc (M, V)
\end{CD}
\]
\\
The second crucial observation is that there is an evident isomorphism between the additive categories $\mathcal{C} (M,V)$ and $\mathcal{C} (U, U \cap V)$.  This isomorphism induces an equivalence in $K$-theory and proves the theorem.

\SSecRef{Coarsely saturated coverings}{CSCOPL}

In order to generalize the excision theorems to fibred categories further in the paper, we need to develop constructions related to coverings of $Y$ and reformulate the excision results.

\begin{DefRef}{HUIKPLO}
Two subsets $A$, $B$ in a proper metric space $X$ are called \textit{coarsely equivalent} if there are numbers $d_{A,B}$, $d_{B,A}$ with
$A \subset B[d_{A,B}]$ and $B \subset A[d_{B,A}]$.
It is clear this is an equivalence relation among subsets.
We will use notation $A \, \| \, B$ for this equivalence.

A family of subsets $\mathcal{A}$ is called \textit{coarsely saturated} if it is maximal with respect to this equivalence relation.
Given a subset $A$, let $\mathcal{S}(A)$ be the smallest boundedly saturated family containing $A$.

If $\mathcal{A}$ is a coarsely saturated family, define $\Knc (\mathcal{A})$ to be
\[
\hocolim{A \in \mathcal{A}} \Knc (A).
\]
\end{DefRef}

The following is immediate from functoriality.

\begin{PropRef}{FGTY}
If $A$, $B$ are coarsely equivalent subsets then $\mathcal{C} (A)$, $\mathcal{C} (B)$ are equivalent categories, and so $\Knc (A)$, $\Knc (B)$ are weakly equivalent spectra.
For all subsets $A$,
$\mathcal{C} (X)_{<A}$ and $\mathcal{C} (A)$ are equivalent categories,
and
$\Knc (X)_{<A} \simeq \Knc (\mathcal{S}(A)) \simeq \Knc (A)$.
\end{PropRef}

\begin{DefRef}{HYUEWS}
A collection of subsets $\mathcal{U} = \{ U_i \}$ is a \textit{coarse covering} of $X$ if $X = \bigcup S_i$ for some $S_i \in \mathcal{S}(U_i)$.
Similarly, $\mathcal{U} = \{ \mathcal{A}_i \}$ is a \textit{coarse covering} by coarsely saturated families if for some (and therefore any) choice of subsets $A_i \in \mathcal{A}_i$, $\{ A_i \}$ is a coarse covering in the sense above.

Recall that a pair of subsets $A$, $B$ in a proper metric space $X$ are \textit{coarsely antithetic} if
for any two numbers $d_1$ and $d_2$ there is a third number $d$ such that
\[
A[d_1] \cap B[d_2] \subset (A \cap B) [d].
\]
We will write $A \, \natural \, B$ to indicate that $A$ and $B$ are coarsely antithetic.

Given two subsets $A$ and $B$, define
\[
\mathcal{S}(A,B) = \{ A' \cap B' \, \vert \, A' \in \mathcal{S}(A), B' \in \mathcal{S}(B), A' \, \natural \, B' \}.
\]
\end{DefRef}

\begin{PropRef}{UYNBVC}
$\mathcal{S}(A,B)$ is a coarsely saturated family.
\end{PropRef}

\begin{proof}
Suppose $A_1$, $A'_1$ and $A_2$, $A'_2$ are two coarsely antithetic pairs, and
$A_1 \subset A_2 [d_{12}]$, $A'_1 \subset A'_2 [d'_{12}]$ for some $d_{12}$ and $d'_{12}$.
Then
\[
A_1 \cap A'_1 \subset A_2 [d_{12}] \cap A'_2 [d'_{12}] \subset (A_2 \cap A'_2) [d]
\]
for some $d$.
\end{proof}

\begin{PropRef}{FGTY2}
If $U$ and $T$ are coarsely antithetic then
\[
K (X)_{<U,T} \simeq K( \mathcal{S}(U \cap T)) \simeq K(U \cap T).
\]
\end{PropRef}

There is the obvious generalization of the constructions and propositions to the case of a finite number of subsets of $X$.

\begin{DefRef}{KIOBXSDZ}
We write $A_1 \, \natural \, A_2 \, \natural \, \ldots \, \natural \, A_k$ if for arbitrary $d_i$ there is a number $d$ so that
\[
A_1 [d_1] \cap A_2 [d_2] \cap \ldots \cap A_k [d_k] \subset (A_1 \cap A_2 \cap \ldots \cap A_k) [d]
\]
and define
\[
\mathcal{S}(A_1, A_2, \ldots , A_k) = \{ A'_1 \cap A'_2 \cap \ldots \cap A'_k \, \vert \, A'_i \in \mathcal{S}(A_i), A_1 \, \natural \, A_2 \, \natural \, \ldots \, \natural \, A_k \}.
\]
Equivalently, identifying any coarsely saturated family $\mathcal{A}$ with $\mathcal{S}(A)$ for $A \in \mathcal{A}$, one has
$\mathcal{S}(\mathcal{A}_1, \mathcal{A}_2, \ldots , \mathcal{A}_k)$.

We will refer to $\mathcal{S}(\mathcal{A}_1, \mathcal{A}_2, \ldots , \mathcal{A}_k)$ as the coarse intersection of $\mathcal{A}_1, \mathcal{A}_2, \ldots , \mathcal{A}_k$.
A coarse covering $\mathcal{U}$ is \textit{closed under coarse intersections} if all coarse intersections $\mathcal{S}(\mathcal{A}_1, \mathcal{A}_2, \ldots , \mathcal{A}_k)$ are nonempty and are contained in $\mathcal{U}$.
If $\mathcal{U}$ is a given coarse covering, the smallest coarse covering that is closed under coarse intersections and contains
$\mathcal{U}$ will be called the \textit{closure} of $\mathcal{U}$ under coarse intersections.
\end{DefRef}

The inclusions induce the diagrams of spectra $\{ K (Y)_{<A} \}$ and
$\{ K^{\Gamma}_i (Y,Y')^{\Gamma}_{<A} \}$ for representatives $A$ in $\mathcal{A} \in \mathcal{U}$.

\begin{DefRef}{HYUN}
Suppose $\mathcal{U}$ is a coarse covering of $Y$ closed under coarse intersections.
We define the homotopy pushouts
\[
\mathcal{K} (Y; \mathcal{U}) = \hocolim{\mathcal{A} \in \mathcal{U}} K(Y)_{<A}.
\]
\end{DefRef}

The following result is equivalent to the Bounded Excision Theorem \refT{BddExc}.

\begin{ThmRef}{PPPOI}
If $\mathcal{U}$ is a finite coarse covering of $Y$ closed under coarse intersections, then there is a weak equivalence
\[
\mathcal{K} (Y; \mathcal{U}) \simeq K(Y).
\]
\end{ThmRef}

\begin{proof}
Apply the Theorem \refT{BddExc} inductively to the sets in $\mathcal{U}$.
\end{proof}

\SSecRef{Some coarse geometry and functoriality}{SCGTZA}

When we develop the equivariant theory, we will want to consider group actions by maps that are more general than isometries.
Let $X$ and $Y$ be proper metric spaces with metric functions $d_X$ and $d_Y$.

\begin{DefRef}{Bornol2}
A map $f \colon X \to Y$
of proper metric spaces is \textit{uniformly expansive} if there is a real positive function $l$
such that
\begin{equation}
d_X (x_1, x_2) \le r \ \Longrightarrow
d_Y (f(x_1), f(x_2)) \le l(r).  \notag
\end{equation}

This is the same concept as \textit{uniformly continuous} maps in Pedersen-Weibel \cite{ePcW:85} or \textit{bornologous} maps in Roe \cite{jR:03}.

A map $f \colon X \to Y$
of proper metric spaces is \textit{proper} if $f^{-1} (S)$ is a bounded subset of $X$ for
each bounded subset $S$ of $Y$.
We say $f$ is a \textit{coarse map} if it is proper and uniformly expansive.

\begin{ThmRef}{JARI}
Coarse maps between proper metric spaces induce additive functors between bounded categories.  Bounded maps induce additive equivalences.
\end{ThmRef}

\begin{proof}
Let $f \colon X \to Y$ be coarse.
We want to induce an additive functor $f_{\ast} \colon \mathcal{C} (X,R) \to \mathcal{C} (Y,R)$.
On objects, the functor is induced by the assignment
\[
(f_{\ast} F)_y = \bigoplus_{x \in f^{-1} (y)} F_x.
\]
Since $f$ is proper, $f^{-1} (y)$ is a bounded set for all $y$ in $Y$.
So the direct sum in the formula is finite, and $(f_{\ast} F)_y$ is a finitely generated free $R$-module.
If $S \subset Y$ is a bounded subset then $f^{-1} (S)$ is bounded.
There are finitely many $F_z \ne 0$  for $z \in f^{-1} (S)$ and therefore finitely many $(f_{\ast} F)_y \ne 0$ for $y \in S$.
This shows $f_{\ast} F$ is locally finite.

Notice that
\[
f_{\ast} F = \bigoplus_{y \in Y} (f_{\ast} F)_y = \bigoplus_{y \in Y} \bigoplus_{z \in f^{-1} (y)} F_z = F.
\]

Suppose we are given a morphism $\phi \colon F \to G$ in $\mathcal{C} (X,R)$.
Interpreting $f_{\ast} F$ and $f_{\ast} G$ as the same $R$-modules, as in the formula above, we define $f_{\ast} \phi \colon f_{\ast} F \to f_{\ast} G$ equal to $\phi$.
We must check that $f_{\ast} \phi$ is bounded.
Suppose $\phi$ is bounded by $D$, and $f$ is $l$-coarse.  We claim that $f_{\ast} \phi$ is bounded by $l (D)$.  Indeed, if $d_{Y} (y,y') > l (D)$ then $d_X (x,x') > D$ for all $x$, $x' \in X$ such that $f(x) = y$ and $f(x') = y'$.  So all components $\phi_{x,x'} = 0$, therefore all components $(f_{\ast} \phi)_{y,y'} = 0$.
\end{proof}

\begin{CorRef}{HMLOWQ}
$\Knc$ is a functor from the category of proper metric spaces and coarse maps to the category of spectra.
\end{CorRef}

The map $f$ is a \textit{coarse equivalence} if there is a coarse map $g \colon Y \to X$ such
that $f \circ g$ and $g \circ f$ are bounded maps.
\end{DefRef}

\begin{Exs} \label{POIU}
Any \textit{bounded} function $f \colon X \to X$, with
$d_X (x, f(x)) \le D$ for all $x \in X$ and a fixed $D \ge 0$, is coarse.
In fact, it is a coarse
equivalence using $l (r) = r + 2D$ for both $f$ and its coarse inverse.

The isometric embedding of a metric subspace is a coarse map.
An isometry, which is a bijective isometric map, is a coarse equivalence.
An isometric embedding onto a subspace that has the property that its bounded enlargement is the whole target metric space is also a coarse equivalence.

A quasi-isometry $f \colon X \to Y$ onto a subset $U \subset Y$ such that for some number $s \ge 0$ we have $U[s]=Y$ is a coarse equivalence \cite{jR:03}.
\end{Exs}

The following definition makes precise a useful class of metrics one has on a finitely generated group.

\begin{DefRef}{WORD}
The \textit{word-length metric} $d = d_{\Omega}$ on a group
$\Gamma$ with a fixed finite generating set $\Omega$ closed under taking inverses is the length metric
induced from the condition that $d (\gamma, \gamma \omega) =1$, whenever $\gamma \in \Gamma$ and $\omega \in \Omega$.  In other words, $d (\alpha, \beta)$ is the minimal length $t$ of sequences $\alpha = \gamma_0, \gamma_1, \ldots  \gamma_t = \beta$ in $\Gamma$ where each consecutive pair of elements differs by right multiplication by an $\omega$ from $\Omega$.
This metric makes $\Gamma$ a proper metric space with a free action by $\Gamma$ via left
multiplication.
\end{DefRef}

If one considers a different choice of a finite generating set $\Omega'$, it is well-known that the identity map on the group with the two metrics $d_{\Omega}$ and $d_{\Omega'}$ is a quasi-isometry.  We see from the combination of Corollary \refC{HMLOWQ} and the examples above that the bounded $K$-theory of $\Gamma$ is independent from the choice of a finite generating set, up to an equivalence.

The following fact is known as ``Milnor's lemma''.

\begin{ThmRefName}{ML}{Shvarts, Milnor}
Suppose $X$ is a length metric space and $\Gamma$ is a group acting properly and cocompactly by isometries on $X$.  Then $\Gamma$ is coarsely equivalent to $X$.
\end{ThmRefName}

This is the situation, for example, when $X$ is a compact smooth Riemannian manifold with the usual action by its fundamental group $\Gamma$.
The coarse equivalence is given by the map $\gamma \mapsto \gamma x_0$ for any fixed base point $x_0$ of $X$.

\SecRef{$K$-theory with fibred control}{KTFCMNB}

\SSecRef{Definition of fibred control}{Basic}

When $M$ is the product of two proper metric spaces $X \times Y$ with the metric
\[
  d_{\textrm{max}} ((x_1,y_1),(x_2,y_2)) = \max \{ d_X (x_1,x_2),
  d_Y (y_1,y_2) \},
\]
one has the bounded category of geometric $R$-modules $\mathcal{C} (X \times Y, R)$ of Pedersen--Weibel as in Definition \refD{PWcats} with the associated $K$-theory spectrum $K (X \times Y, R)$.
We construct another category associated to the pair $(X, Y)$.

\begin{DefRef}{FBCcat}
The new category $\mathcal{C}_X (Y)$ has the same objects as $\mathcal{C} (X \times Y, R)$ but a weaker control condition on the morphisms.
For any function $f \colon [0, +\infty) \to [0, +\infty)$ and a real number $D \ge 0$, define
\[
N (D,f) (x,y) = x[D] \times y[f(d(x,x_0))],
\]
the $(D,f)$-neighborhood of $(x,y)$ in $X \times Y$.
A homomorphism $\phi \colon F \to G$ is called $(D,f)$-\textit{bounded} if the components $F_{(x,y)} \to G_{(x',y')}$ are zero maps for $(x',y')$ outside of the $(D,f)$-neighborhood of $(x,y)$.  These are the morphisms of $\mathcal{C}_X (Y)$.
It is easy to see that the definition of the category $\mathcal{C}_X (Y)$ is independent of the choice of $x_0$ in $X$.
We will refer to this category of geometric modules over the product as a \textit{fibred bounded category}.
\end{DefRef}

The connective $K$-theory spectrum $K_X (Y)$ is the spectrum associated to the isomorphism category of $\mathcal{C}_X (Y)$.

\begin{RemRef}{FtF}
It follows from Example 1.2.2 of Pedersen-Weibel \cite{ePcW:85} that in general the fibred bounded category $\mathcal{C}_X (Y)$ is not
isomorphic to $\mathcal{C} (X \times Y, R)$.
The proper generality of that work, as explained in \cite{ePcW:89}, has a general additive category $\mathcal{A}$ embedded in a cocomplete additive category, generalizing the setting of free finitely generated $R$-modules as a subcategory of all free $R$-modules.
All of the excision results of Pedersen-Weibel hold for $\mathcal{C} (X, \mathcal{A})$.
In these terms, the category $\mathcal{C}_X (Y)$ is isomorphic to the category $\mathcal{C} (X, \mathcal{A})$, where $\mathcal{A} = \mathcal{C} (Y, R)$.

The difference between $\mathcal{C}_X (Y)$ and $\mathcal{C} (X \times Y, R)$ is made to disappear in \cite{ePcW:85}
by making $\mathcal{C} (Y,R)$ ``remember the filtration'' of morphisms when viewed
as a filtered additive category with
$\Hom_D (F,G) = \{ \phi \in \Hom (F,G) \mid \fil (\phi) \le D \}$.
Identifying a small category with its set of morphisms, one can
think of the bounded category as
\[
  \mathcal{C} (Y,R) = \colim{D \in \mathbb{R}} \mathcal{C}_D (Y,R),
\]
where $\mathcal{C}_D (Y,R) = \Hom_D (\mathcal{C} (Y,R))$ is the
collection of all $\Hom_D (F,G)$. Now we have
\[
\mathcal{C} (X \times Y, R)
= \colim{D \in \mathbb{R}} \mathcal{C} (X,
\mathcal{C}_D (Y,R)).
\]

There is an exact embedding $\iota \colon \mathcal{C} (X \times Y, R) \to \mathcal{C}_X (Y)$ inducing the map of $K$-theory spectra $K(\iota) \colon K (X \times Y, R) \to K_X (Y)$.
\end{RemRef}

We want to develop some results for a variety of categories with fibred bounded control where the Karoubi filtration techniques suffice.

\begin{NotRef}{JUODQ}
Let
\begin{align}
\mathcal{C}_k = \ &\mathcal{C}_{X} (Y \times \mathbb{R}^k), \notag \\
\mathcal{C}^{+}_k = \ &\mathcal{C}_{X} (Y \times \mathbb{R}^{k-1} \times [0,+\infty)), \notag \\
\mathcal{C}^{-}_k = \ &\mathcal{C}_{X} (Y \times \mathbb{R}^{k-1} \times (-\infty,0]). \notag
\end{align}
We will also use the notation
\begin{align}
\mathcal{C}^{<+}_k =& \colim{D \ge 0} \mathcal{C}_{X} (Y \times \mathbb{R}^{k-1} \times [-D,+\infty)), \notag \\
\mathcal{C}^{<-}_k =& \colim{D \ge 0} \mathcal{C}_{X} (Y \times \mathbb{R}^{k-1} \times (-\infty,D]), \notag \\
\mathcal{C}^{<0}_k =& \colim{D \ge 0} \mathcal{C}_{X} (Y \times \mathbb{R}^{k-1} \times [-D,D]). \notag
\end{align}
\end{NotRef}
Clearly $\mathcal{C}_k$ is $\mathcal{C}^{<-}_k$-filtered and that $\mathcal{C}^{<+}_k$ is $\mathcal{C}^{<0}_k$-filtered.  There are isomorphisms $\mathcal{C}^{<0}_k \cong \mathcal{C}_{k-1}$, $\mathcal{C}^{<-}_k \cong \mathcal{C}^{-}_k$, and $\mathcal{C}_k/ \mathcal{C}^{<-}_k \cong \mathcal{C}^{<+}_k/ \mathcal{C}^{<0}_k$.  By Theorem \refT{CPF}, the commutative diagram
\[
\begin{CD}
K ((\mathcal{C}^{<0}_k)^{\wedge K}) @>>> K (\mathcal{C}^{<+}_k) @>>> K (\mathcal{C}^{<+}_k/ \mathcal{C}^{<0}_k) \\
@VVV @VVV @VV{\cong}V \\
K ((\mathcal{C}^{<-}_k)^{\wedge K}) @>>> K (\mathcal{C}_k) @>>> K (\mathcal{C}_k/ \mathcal{C}^{<-}_k)
\end{CD}
\]
where all maps are induced by inclusions on objects, is in fact a map of homotopy fibrations.  The categories $\mathcal{C}^{<+}_k$ and $\mathcal{C}^{<-}_k$ are flasque, that is, possess an endofunctor $\Sh$ such that $\Sh (F) \cong F \oplus \Sh (F)$, which can be seen by the usual Eilenberg swindle argument.
Therefore $K (\mathcal{C}^{<+}_k)$ and $K (\mathcal{C}^{<-}_k)$ are contractible by the
Additivity Theorem, cf.\ Pedersen--Weibel \cite{ePcW:85}.   This gives a map $K (\mathcal{C}_{k-1}) \to \Omega  K (\mathcal{C}_k)$ which induces isomorphisms of $K$-groups in positive dimensions.

\begin{DefRef}{POIYSNHU}
We define the \textit{nonconnective fibred bounded $K$-theory} over the pair $(X,Y)$ as the spectrum
\[
\KncX (Y) \overset{ \text{def} }{=} \hocolim{k>0}
\Omega^{k} K (\mathcal{C}_k).
\]
\end{DefRef}

If $Y$ is the single point space then the delooping $\KncX (\point)$ is clearly equivalent to
the nonconnective delooping $\Knc (X,R)$ of Pedersen--Weibel, reviewed in Theorem \refT{Deloop}, via the map
$K(\iota) \colon \Knc (X \times \point, R) \rightarrow \KncX (\point)$.

\SSecRef{Fibrewise coarsely saturated coverings}{CSFCOPL}

Following the blueprint from section \refS{KFPWC}, we now describe proper systems of coverings and prove the fibred excision theorem for these systems.

Suppose we are given a finite coarse antithetic covering $\mathcal{U}$ of $Y$ closed under coarse intersections and of cardinality $s$.
The coarsely saturated families $\mathcal{A}_i$ which are members of $\mathcal{U}$ are partially ordered by inclusion.
In fact, the union of the families $\mathcal{A}_i$ forms the set $\mathcal{A}$ closed under intersections.

We will use the notation $(X,Y)$ for the product $X \times Y$ when we want to distinguish the roles the factors play.

\begin{DefRef}{HUIKPLO2}
Two subsets $A$, $B$ of $(X,Y)$ are called \textit{coarsely equivalent} if there is a set of enlargement data $(K,k)$ such that
$A \subset B[K,k_{x_0}]$ and $B \subset A[K,k_{x_0}]$.
It is an equivalence relation among subsets.
We will again use notation $A \, \| \, B$ for this equivalence generalizing of the notation from
Definition \refD{HUIKPLO}.

A family of subsets $\mathcal{A}$ is called \textit{coarsely saturated} if it is maximal with respect to this equivalence relation.
Given a subset $A$, let $\mathcal{S}(A)$ be the smallest boundedly saturated family containing $A$.

A collection of subsets $\mathcal{U} = \{ U_i \}$ is a \textit{coarse covering} of $(X,Y)$ if $(X,Y) = \bigcup S_i$ for some $S_i \in \mathcal{S}(U_i)$.
Similarly, $\mathcal{U} = \{ \mathcal{A}_i \}$ is a \textit{coarse covering} by coarsely saturated families if for some (and therefore any) choice of subsets $A_i \in \mathcal{A}_i$, $\{ A_i \}$ is a coarse covering in the sense above.

A pair of subsets $A$, $B$ of $(X,Y)$ is called  \textit{coarsely antithetic} if
for any two sets of enlargement data $(D_1, d_1)$ and $(D_2, d_2)$ there is a third set $(D, d)$ such that
\[
A[D_1, (d_1)_{x_0}] \cap B[D_2, (d_2)_{x_0}] \subset (A \cap B) [D, d_{x_0}].
\]
We will write $A \, \natural \, B$ to indicate that $A$ and $B$ are coarsely antithetic.

Given two subsets $A$ and $B$, we define
\[
\mathcal{S}(A,B) = \{ A' \cap B' \, \vert \, A' \in \mathcal{S}(A), B' \in \mathcal{S}(B), A' \, \natural \, B' \}.
\]
It is easy to see that $\mathcal{S}(A,B)$ is a coarsely saturated family, cf.~Proposition \refP{UYNBVC}.

There is the straightforward generalization to the case of a finite number of subsets of $(X,Y)$.
Again, we write $A_1 \, \natural \, A_2 \, \natural \, \ldots \, \natural \, A_k$ if for arbitrary sets of data $(D_i,d_i)$ there is a set of enlargement data $(D,d)$ so that
\begin{multline*}
A_1 [D_1,(d_1)_{x_0}] \cap A_2 [D_2,(d_2)_{x_0}] \cap \ldots \cap A_k [D_k,(d_k)_{x_0}] \\
\subset (A_1 \cap A_2 \cap \ldots \cap A_k) [D,d_{x_0}]
\end{multline*}
and define
\[
\mathcal{S}(A_1, A_2, \ldots , A_k) = \{ A'_1 \cap A'_2 \cap \ldots \cap A'_k \, \vert \, A'_i \in \mathcal{S}(A_i), A_1 \, \natural \, A_2 \, \natural \, \ldots \, \natural \, A_k \}.
\]
Identifying any coarsely saturated family $\mathcal{A}$ with $\mathcal{S}(A)$ for $A \in \mathcal{A}$, one has
the coarse saturated family
$\mathcal{S}(\mathcal{A}_1, \mathcal{A}_2, \ldots , \mathcal{A}_k)$.
We will refer to $\mathcal{S}(\mathcal{A}_1, \mathcal{A}_2, \ldots , \mathcal{A}_k)$ as the
\textit{coarse intersection} of $\mathcal{A}_1, \mathcal{A}_2, \ldots , \mathcal{A}_k$.
A coarse covering $\mathcal{U}$ is \textit{closed under coarse intersections} if all coarse intersections $\mathcal{S}(\mathcal{A}_1, \mathcal{A}_2, \ldots , \mathcal{A}_k)$ are nonempty and are contained in $\mathcal{U}$.
\end{DefRef}

\begin{PropRef}{UIVCXZSAW}
If $\mathcal{U}$ is a coarse antithetic covering of $Y$ then $(X,\mathcal{U})$ consisting of subsets $(X,U)$, $U \in \mathcal{U}$, is a coarse antithetic covering of $(X,Y)$.  If $\mathcal{U}$ is closed under coarse intersections, $(X,\mathcal{U})$ is closed under coarse intersections.
\end{PropRef}

\begin{proof}
Left to the reader.
\end{proof}

\SSecRef{Fibrewise excision theorems}{FETDSLK}

Bounded excision theorems of section \refSS{BETRGCV} can be adapted to the fibred setting with 
$\KncX (Y)$.
Suppose $Y_1$
and $Y_2$ are mutually antithetic subsets of a proper metric space $Y$, and $Y = Y_1
\cup Y_2$.
Consider the coarse covering $\mathcal{U}$ of $Y$ by $\mathcal{S}(Y_1)$, $\mathcal{S}(Y_2)$, and $\mathcal{S}(Y_1,Y_2)$.

\begin{NotRef}{ZXCVRE}
	Let us use the general notation $\mathcal{C}_X
(Y)_{<C}$, for a subset $C$ of $Y$, to denote the full subcategory of $\mathcal{C}_X (Y)$ on the objects with supports $Z \subset X \times Y$ which are coarsely equivalent to the subset $X \times C$.
	
	We will also use the shorthand notation $\mathcal{C} = \mathcal{C}_X (Y)_{<\mathcal{U}}$, $\mathcal{C}_i= \mathcal{C}_X
(Y)_{<Y_i}$ for $i=1$ or $2$, and $\mathcal{C}_{12}$ for the intersection
$\mathcal{C}_1 \cap \mathcal{C}_2$.
\end{NotRef}

There is a commutative diagram
\[
\begin{CD}
K (\mathcal{C}_{12}) @>>> K (\mathcal{C}_1) @>>> K ({\mathcal{C}_1}/{\mathcal{C}_{12}}) \\
@VVV @VVV @VV{K(I)}V \\
K (\mathcal{C}_2) @>>> K (\mathcal{C}) @>>> K ({\mathcal{C}}/{\mathcal{C}_2})
\end{CD} \tag{$\natural$}
\]
where the rows are homotopy fibrations from Theorem \refT{CPF} and $I \colon {\mathcal{C}_1}/{\mathcal{C}_{12}} \to {\mathcal{C}}/{\mathcal{C}_2}$ is the functor induced from the exact inclusion $I \colon \mathcal{C}_1 \to \mathcal{C}$.
It clear that $I$ is again an additive isomorphism.

The subcategory $\mathcal{C}_{X} (Y \times \mathbb{R}^k)_{<C \times \mathbb{R}^k}$ is evidently a Karoubi subcategory of $\mathcal{C}_{X} (Y \times \mathbb{R}^k)$ for any choice of the subset $C \subset Y$.
We define
\[
\KncX (Y)_{<C} \overset{ \text{def} }{=} \hocolim{k>0}
\Omega^{k} {K}_{X} (Y \times \mathbb{R}^k)_{<C \times \mathbb{R}^k}.
\]
Using the methods above, one easily obtains the weak equivalence
\begin{gather*}
\KncX (Y)_{<C} \simeq \KncX (C).
\end{gather*}
We also define
\[
\KncX (Y)_{<Y_1, Y_2} \overset{ \text{def} }{=} \hocolim{k>0}
\Omega^{k} {K}_{X} (Y \times \mathbb{R}^k)_{<Y_1 \times \mathbb{R}^k, \, Y_2 \times \mathbb{R}^k}.
\]

We are now able to follow the proof of Theorem \refT{BddExc} to get the following.

\begin{ThmRefName}{Exc2}{Fibrewise bounded excision}
Suppose $Y_1$ and $Y_2$ are subsets of a metric space $Y$, and $Y = Y_1 \cup Y_2$.
There is a homotopy pushout diagram of spectra
\[
\begin{CD}
\KncX (Y)_{<Y_1,Y_2} @>>> \KncX (Y)_{<Y_1} \\
@VVV @VVV \\
\KncX (Y)_{<Y_2} @>>> \KncX (Y)
\end{CD}
\]
where the maps are induced from the exact inclusions.
If $Y_1$ and $Y_2$ are mutually antithetic subsets of $Y$,
there is a homotopy pushout
\[
\begin{CD}
\KncX (Y_1 \cap Y_2) @>>> \KncX (Y_1) \\
@VVV @VVV \\
\KncX (Y_2) @>>> \KncX (Y)
\end{CD}
\]
\end{ThmRefName}

Suppose we are given a finite coarse antithetic covering $\mathcal{U}$ of $Y$ closed under coarse intersections and of cardinality $s$.
The coarsely saturated families $\mathcal{A}_i$ which are members of $\mathcal{U}$ are partially ordered by inclusion.
In fact, the union of the families $\mathcal{A}_i$ forms the set $\mathcal{A}$ closed under intersections.

\begin{PropRef}{UIVCXZSAW}
If $\mathcal{U}$ is a coarse antithetic covering of $Y$ then $(X,\mathcal{U})$ consisting of subsets $(X,U)$, $U \in \mathcal{U}$, is a coarse antithetic covering of $(X,Y)$.  If $\mathcal{U}$ is closed under coarse intersections, $(X,\mathcal{U})$ is closed under coarse intersections.
\end{PropRef}

\begin{proof}
Left to the reader.
\end{proof}

\begin{CorRef}{gnaloi}
Suppose $\mathcal{U}$ is a finite coarse covering of $Y$ closed under coarse intersections.
We can define the homotopy pushout
\[
\mathcal{K}_X (Y; \mathcal{U}) = \hocolim{A \in \mathcal{A} \in \mathcal{U}} \KncX (Y)_{<A}.
\]
Then there is a weak equivalence
$\mathcal{K}_X (Y; \mathcal{U}) \simeq \KncX (Y)$.
\end{CorRef}

\begin{proof}
Apply Theorem \refT{Exc2} inductively to the sets in $\mathcal{U}$.
\end{proof}

There is a relative version of fibred $K$-theory and the corresponding Fiberwise Excision Theorem.

\begin{DefRef}{GEPCWpr}
Let $Y' \in \mathcal{A}$ as part of a coarse covering $\mathcal{U}$ of $Y$.  Let $\mathcal{C} = \mathcal{C}_X (Y)_{< \mathcal{U}}$ and $\mathcal{Y}' = \mathcal{C}_X (Y)_{<Y'}$.  The category $\mathcal{C}_X (Y,Y')$ is the quotient category $\mathcal{C}/\mathcal{Y}'$.

It is now straightforward to define
\[
\KncX (Y,Y') {=} \hocolim{k>0}
\Omega^{k} {K}_{X} (Y \times \mathbb{R}^k, Y' \times \mathbb{R}^k),
\]
\[
\KncX (Y,Y')_{<C} {=} \hocolim{k>0}
\Omega^{k} {K}_{X} (Y \times \mathbb{R}^k, Y' \times \mathbb{R}^k)_{<C \times \mathbb{R}^k},
\]
and
\[
\KncX (Y,Y')_{<C_1, C_2} {=} \hocolim{k>0}
\Omega^{k} {K}_{X} (Y \times \mathbb{R}^k, Y' \times \mathbb{R}^k)_{<C_1 \times \mathbb{R}^k, \, C_2 \times \mathbb{R}^k}.
\]
\end{DefRef}

The theory developed in this section is easily relativized to give the following excision theorems.

\begin{ThmRefName}{ExRel}{Relative fibred excision}
If $Y$ is the union of two subsets $Y_1$ and $Y_2$,
there is a homotopy pushout of spectra
\[
\begin{CD}
\KncX (Y,Y')_{<Y_1,Y_2} @>>> \KncX (Y,Y')_{<Y_1} \\
@VVV @VVV \\
\KncX (Y,Y')_{<Y_2} @>>> \KncX (Y,Y')
\end{CD}
\]
where the maps are induced from the exact inclusions.
In fact, if $Y$ is the union of two mutually antithetic subsets $Y_1$ and $Y_2$, and $Y'$ is antithetic to both $Y_1$ and $Y_2$,
there is a homotopy pushout
\[
\begin{CD}
\KncX (Y_1 \cap Y_2, Y_1 \cap Y_2 \cap Y') @>>> \KncX (Y_1, Y_1 \cap Y') \\
@VVV @VVV \\
\KncX (Y_2, Y_2 \cap Y') @>>> \KncX (Y,Y')
\end{CD}
\]
More generally, we can define the homotopy pushout
\[
\mathcal{K}_X (Y,Y'; \mathcal{Y}) = \hocolim{A \in \mathcal{A} \in \mathcal{Y}} \KncX (Y,Y')_{<A}.
\]
Then there is a weak equivalence
\[
\mathcal{K}_X (Y,Y'; \mathcal{Y}) \simeq \KncX (Y,Y').
\]
\end{ThmRefName}

\begin{ThmRef}{ExcExcPrep}
Given a subset $U$ of $Y'$, there is an equivalence
\[
\KncX (Y,Y') \simeq \KncX (Y-U,Y'-U).
\]
\end{ThmRef}

\begin{proof}
Consider the setup of Theorem \refT{Exc2} with $Y_1 = Y-U$ and $Y_2=Y'$, then the map
\[
\frac{\mathcal{C}_X (Y)_{<(Y-U)}}{\mathcal{C}_X (Y)_{<(Y-U)} \cap \mathcal{C}_X (Y)_{<Y'}} \longrightarrow \frac{\mathcal{C}_X (Y)}{\mathcal{C}_X (Y)_{<Y'}}
\]
is an isomorphism and so induces a weak equivalence on the level of $K$-theory.
Notice that, since $U$ is a subset of $Y'$,
\[
{\mathcal{C}_X (Y)_{<(Y-U)} \cap \mathcal{C}_X (Y)_{<Y'}} = \mathcal{C}_X (Y)_{<(Y'-U)}.
\]
Now the maps of quotients
\[
\frac{\mathcal{C}_X (Y)}{\mathcal{C}_X (Y')}  \longrightarrow  \frac{\mathcal{C}_X (Y)}{\mathcal{C}_X (Y)_{<Y'}}
\ \ \,
\mathrm{and}
\ \ \,
  \frac{\mathcal{C}_X (Y-U)}{\mathcal{C}_X (Y'-U)}
  \longrightarrow  
  \frac{\mathcal{C}_X (Y)_{<(Y-U)}}{\mathcal{C}_X (Y)_{<(Y'-U)}}
\]
also induce weak equivalences.
Their composition gives the required equivalence.
\end{proof}

\SecRef{Equivariant theory}{JKASEQT}

\SSecRef{Basic equivariant theory}{OET}

The classical situation is a proper metric space $M$ with a free (left) $\Gamma$-action by isometries.  
Clearly,
there is a natural action of $\Gamma$ on the geometric modules $\mathcal{C} (M)$ and
therefore on $\Knc (M)$.  Formally, this action is also free.
A different equivariant bounded $K$-theory with useful fixed point spectra is constructed as follows.

\begin{DefRef}{EqCatDef}
Let $\EGamma$ be the category with the object set $\Gamma$ and the
unique morphism $\mu \colon \gamma_1 \to \gamma_2$ for any pair
$\gamma_1$, $\gamma_2 \in \Gamma$. There is a left action by $\Gamma$
on $\EGamma$ induced from the left multiplication in $\Gamma$.
\end{DefRef}

If $\mathcal{C}$ is a small category with left $\Gamma$-action,
then the functor category
$\Fun(\EGamma,\mathcal{C})$ is a category with the
left $\Gamma$-action given on objects by
$\gamma(F)(\gamma')=\gamma F (\gamma^{-1} \gamma')$
and
$\gamma(F)(\mu)=\gamma F (\gamma^{-1} \mu)$.
It is always
nonequivariantly equivalent to $\mathcal{C}$.

\begin{DefRef}{LaxLimit}
The subcategory of equivariant functors and
equivariant natural transformations in $\Fun(\EGamma,\mathcal{C})$
is the fixed subcategory $\Fun(\EGamma,\mathcal{C})^{\Gamma}$ known as
the \textit{lax limit of the action}, introduced and called so by Thomason in \cite{rT:83}.

Thomason also gave the following explicit description of the lax limit.
The objects of
$\Fun(\EGamma,\mathcal{C})^{\Gamma}$ can
be thought of as pairs
$(F,\psi)$ where $F \in \mathcal{C}$ and $\psi$ is a function on
$\Gamma$ with $\psi (\gamma) \in \Hom (F,\gamma F)$ subject to the condition 
$\psi(1) = 1$ and the cocycle identity $\psi (\gamma_1 \gamma_2) =
\gamma_1 \psi(\gamma_2)  \circ \psi (\gamma_1)$.
These conditions imply that $\psi (\gamma)$ is always an
isomorphism. 

The set of morphisms $(F,\psi) \to
(F',\psi')$ consists of the morphisms $\phi \colon F \to F'$ in
$\mathcal{C}$ such that the squares
\[
\begin{CD}
F @>{\psi (\gamma)}>> \gamma F \\
@V{\phi}VV @VV{\gamma \phi}V \\
F' @>{\psi' (\gamma)}>> \gamma F'
\end{CD}
\]
commute for all $\gamma \in \Gamma$.
\end{DefRef}

Specializing this definition to the case of $\mathcal{C} = \mathcal{C} (M)$, we will use the notation
$\mathcal{C}^{\Gamma} (M)$ for the equivariant theory $\Fun(\EGamma,\mathcal{C} (M))$.
It turns out that the following adjustment makes the fixed points better behaved.
Notice that
$\mathcal{C} (M)$ contains the family of isomorphisms $\phi$ such that $\phi$ and $\phi^{-1}$ are bounded by $0$.
We will express this property by saying that the filtration of $\phi$ is $0$ and writing $\fil (\phi) =0$.
The full
subcategory of functors $\theta \colon \EGamma \to \mathcal{C} (M)$ such that $\fil \theta
(f) = 0$ for all $f$ is invariant under the $\Gamma$-action.

\begin{DefRef}{PrKeq}
Let $\mathcal{C}^{\Gamma,0} (M)$
be the full subcategory of $\mathcal{C}^{\Gamma} (M)$ on the functors sending all
morphisms of $\EGamma$ to filtration $0$ maps.
Then $\mathcal{C}^{\Gamma,0} (M)^{\Gamma}$  is the full subcategory of $\mathcal{C}^{\Gamma} (M)^{\Gamma}$ on the objects
$(F, \psi)$ with $\fil \psi (\gamma) = 0$ for all $\gamma \in \Gamma$.
\end{DefRef}

We define
$K^{\Gamma,0} (M)$ to be the nonconnective delooping of the $K$-theory of the symmetric monoidal category $\mathcal{C}^{\Gamma,0} (M)$.

It is shown in section VI of \cite{gC:95} that the fixed points of the induced $\Gamma$-action on $K^{\Gamma,0} (M)$ is the nonconnective delooping of the $K$-theory of
$\mathcal{C}^{\Gamma,0} (M)^{\Gamma}$.
One also knows the following.

\begin{ThmRefName}{VVV}{Corollary VI.8 of \cite{gC:95}}
If $M$ is a proper metric space and $\Gamma$ acts on $M$ freely, properly discontinuously, cocompactly by isometries, there are weak equivalences
$K^{\Gamma,0} (M)^{\Gamma} \simeq 
\Knc (M/\Gamma, R[\Gamma]) \simeq \Knc (R[\Gamma])$.
\end{ThmRefName}

This theorem applies in two specific cases of interest to us.
One case is of the
fundamental group $\Gamma$ of a closed aspherical manifold acting on the universal cover by
covering transformations.
The second is of $\Gamma$ acting on itself, as a word-length metric space, by left
multiplication.

\begin{RemRef}{NotSame}
The theory $K^{\Gamma,0} (M)^{\Gamma}$ may very well
differ from $K^{\Gamma} (M)^{\Gamma}$.
According to Theorem \refT{VVV}, the fixed point category
of the correct bounded equivariant theory is, for example, the category of free
$R[\Gamma]$-modules when $M = \Gamma$ with the word metric relative to a chosen generating set. 
However,
$\mathcal{C}^{\Gamma} (\Gamma,R)^{\Gamma}$ will
include the $R$-module with a single stalk $F_{\gamma} = R$ for all $\gamma$ and equipped with trivial
$\Gamma$-action, which is not free.
\end{RemRef}

\SSecRef{Equivariant fibred theories}{EETT}

We will treat the group $\Gamma$ equipped with a finite generating set $\Omega$ as a metric space with a word-length metric, as explained in Definition \refD{WORD}.  It is well-known that varying $\Omega$ only changes $\Gamma$ to a
coarsely equivalent metric space.

\begin{DefRefName}{POI1}{Coarse equivariant theories}
We associate two new equivariant theories on metric spaces with a left $\Gamma$-action, both by
isometries and coarse equivalences.  The theory $K_i^{\Gamma}$ is defined only for metric
spaces $Y$ with actions by isometries, while $K_p^{\Gamma}$ for metric
spaces with coarse actions.
\begin{enumerate}
\item
$k^{\Gamma}_i (Y)$ is defined to be the $K$-theory of
$\mathcal{C}^{\Gamma}_i (Y) = \mathcal{C}^{\Gamma,0}(\Gamma \times Y, R)$,
where $\Gamma$ is regarded as a
word-length metric space with isometric $\Gamma$-action given by left multiplication, and $\Gamma
\times Y$ is given the product metric and the diagonal isometric action.
\item
$k^{\Gamma}_p (Y)$ is defined for any metric space $Y$ equipped with a $\Gamma$-action
by coarse equivalences. It is the $K$-theory spectrum attached to a symmetric
monoidal category $\mathcal{C}^{\Gamma}_p (Y)$ with $\Gamma$-action whose objects are given by functors
\[
\theta \colon \EGamma \longrightarrow \mathcal{C}_{\Gamma} (Y) =
\mathcal{C} (\Gamma, \mathcal{C} (Y,R))
\]
with the additional condition that the morphisms $\theta(f)$ are of filtration zero but only as homomorphisms between $R$-modules parametrized over $\Gamma$.
\end{enumerate}
\end{DefRefName}

Now the nonconnective equivariant $K$-theory spectra $K_i^{\Gamma}$ and $K_p^{\Gamma}$ should be the nonconnective deloopings of $k_i^{\Gamma}$ and $k_p^{\Gamma}$.
For example, if we define
\[
\mathcal{C}_{i,k}^{\Gamma} = \mathcal{C}^{\Gamma,0}(\Gamma \times \mathbb{R}^{k} \times Y, R),
\]
where $\Gamma$ acts on the product $\Gamma \times \mathbb{R}^{k} \times Y$ according to $\gamma  (\gamma', x, y) = (\gamma \gamma', x, \gamma(y))$, and
\[
\mathcal{C}_{i,k}^{\Gamma,+} = \mathcal{C}^{\Gamma,0}(\Gamma \times [0,+\infty) \times Y, R), \ \mathrm{etc.,}
\]
then the delooping construction in Definition \refD{POIYSNHU} can be applied verbatim.
Similarly, one can use the $\Gamma$-action on $\Gamma \times \mathbb{R}^{k}$ given by
$\gamma  (\gamma', x) = (\gamma \gamma', x)$ and define $\mathcal{C}_{p,k}^{\Gamma}$ as the symmetric monoidal category of functors
$\theta \colon \EGamma \rightarrow
\mathcal{C} (\Gamma \times \mathbb{R}^{k}, \mathcal{C} (Y,R))$
such that the morphisms $\theta(f)$ are bounded by $0$ as $R$-linear homomorphisms
over $\Gamma \times \mathbb{R}^{k}$.
There are obvious analogues of the categories $\mathcal{C}_{p,k}^{\Gamma,+}$, etc.
If $\ast$ is either of the subscripts $i$ or $p$, and $Y$ is equipped with actions by $\Gamma$ via respectively isometries or coarse equivalences, we obtain equivariant maps
\[
K (\mathcal{C}^{\Gamma}_{\ast, k-1} (Y)) \longrightarrow \Omega K (\mathcal{C}^{\Gamma}_{\ast,k} (Y)).
\]

\begin{DefRef}{EQIOMN}
Let $\ast$ be either of the subscripts $i$ or $p$.
We define
\[
k^{\,\Gamma}_{\ast,k} (Y) = K (\mathcal{C}^{\Gamma}_{\ast,k} (Y))
\]
and the nonconnective equivariant spectra
\[
K^{\Gamma}_{\ast} (Y) \overset{ \text{def} }{=} \hocolim{k>0}
\Omega^{k} k^{\,\Gamma}_{\ast,k} (Y).
\]
The same construction gives for the fixed points
\[
K^{\Gamma}_{\ast} (Y)^{\Gamma} = \hocolim{k>0}
\Omega^{k} k^{\,\Gamma}_{\ast,k} (Y)^{\Gamma}.
\]

Given left $\Gamma$-actions by isometries on metric spaces $X$ and $Y$, there are evident diagonal actions induced on the categories $\mathcal{C} (X \times Y, R)$ and $\mathcal{C}
(X, \mathcal{C} (Y,R))$. The equivariant embedding induces
the equivariant functor
\[
i^{\Gamma} \colon \mathcal{C}^{\Gamma} (X \times Y, R) \longrightarrow
\mathcal{C}^{\Gamma}
(X, \mathcal{C} (Y,R)).
\]
In this case
there is a natural transformation $K^{\Gamma}_i (Y) \to K^{\Gamma}_p (Y)$.
\end{DefRef}

One basic relation between the two equivariant fiberwise theories is through the observation that
in both cases, when $Y$ is a single point space, $\mathcal{C}^{\Gamma}_i (\point)$  and $\mathcal{C}^{\Gamma}_p (\point)$ can be identified
with $\mathcal{C}^{\Gamma,0} (\Gamma, R)$.

Another fact is a special instructive property of the theory $K_p^{\Gamma}$.

\begin{ThmRef}{GpDef}
Suppose $\Gamma$ acts on a metric space $Y$ by bounded coarse equivalences.  Let $Y_0$ be the same metric space but with $\Gamma$ acting trivially by the identity.  Then there is a weak equivalence
\[
\zeta \colon
K_p^{\Gamma} (Y)^{\Gamma} \xrightarrow{ \ \simeq \ } K_p^{\Gamma} (Y_0)^{\Gamma}.
\]
\end{ThmRef}

\begin{proof}
The category $\mathcal{C}_p^{\Gamma} (Y)$ has the left action by $\Gamma$ induced from the diagonal action on $\Gamma \times Y$.
Recall that an object of $\mathcal{C}_p^{\Gamma} (Y)^{\Gamma}$ is determined by an object
$F$ of $\mathcal{C}_{\Gamma} (Y)$ and isomorphisms $\psi (\gamma) \colon F \to {\gamma} F$ which are of filtration $0$ when projected to $\Gamma$.
Given two objects $(F, \{\psi (\gamma) \})$ and $(G, \{\phi (\gamma) \})$, a morphism $\lambda \colon (F, \{\psi (\gamma) \}) \to (G, \{\phi (\gamma) \})$ is given by a morphism $\lambda \colon F \to G$ in $\mathcal{C}_{\Gamma} (Y)$ such that the collection of morphisms $\gamma \lambda \colon {\gamma} F \to {\gamma} G$ satisfies
\[
\psi (\gamma) \circ \lambda =  \gamma \lambda \circ \phi (\gamma)
\]
for all $\gamma$ in $\Gamma$.
Given $(F, \{\psi (\gamma) \})$, define $(F_0, \{\psi_0 (\gamma) \})$ by $F_0 = F$ and $\psi_0 (\gamma)= \id_{F}$ for all $\gamma \in \Gamma$.  Then $\psi (\gamma)^{-1}$ give a natural isomorphism $Z_F$ from $F$ to $F_0$ and induce
an equivalence
\[
\zeta \colon K_p^{\Gamma} (Y)^{\Gamma} \, \simeq \,  K_p^{\Gamma} (Y_0)^{\Gamma}.
\]
Of course, the bound for the isomorphism $\psi (\gamma)^{-1}$ can vary with $\gamma$.
\end{proof}

Notice that this equivalence exists only in the theory $K_p^{\Gamma}$ and not in $K_i^{\Gamma}$ because the homomorphism $\lambda$ we used in the proof is not a bounded homomorphism and is a morphism only in $\mathcal{C}_p^{\Gamma}$.

\SSecRef{Fibrewise equivariant localization and excision}{FLEKth}

Throughout this section we will fix one left action of $\Gamma$ on $Y$ by bounded coarse equivalences.  This means that for each element $\gamma$ the corresponding self-equivalence of $Y$ is a bounded map, but the bound is allowed to vary with $\gamma$.  

We will state results only for the equivariant theory $K_p^{\Gamma}$ where they are the most useful. Some but not all of the statements are also true in $K_i^{\Gamma}$.  However, it is easy to see that all of the statements we make are true for $K_i^{\Gamma}$ if the action is the trivial action by the identity for all $\gamma$.

\begin{DefRef}{KQBdd}
Using the conventions from Notation \refN{ZXCVRE}, let
$\mathcal{C}^{\Gamma}_{p} (Y)_{<Y'}$ be the full subcategory of $\mathcal{C}^{\Gamma}_{p} (Y)$ on objects $\theta$ such that the support of each $\theta (\gamma)$ is contained in a subset coarsely equivalent to $\Gamma \times Y'$.
It is evident that the subcategories give Karoubi filtrations and therefore Karoubi quotients
$\mathcal{C}^{\Gamma}_{p} (Y,Y')$.
The bounded actions of $\Gamma$ extend to the quotients in each case.
Taking $K$-theory of the equivariant symmetric monoidal categories gives the $\Gamma$-equivariant spectra
$k^{\Gamma}_{p} (Y,Y')$.

One can now construct the parametrized versions of the relative module categories
$\mathcal{C}^{\Gamma}_{p,k} (Y,Y')$,
their $K$-theory spectra $k^{\Gamma}_{p,k} (Y,Y')$, and the resulting deloopings.
Thus we obtain the nonconnective $\Gamma$-equivarint spectra
\[
K^{\Gamma}_{p} (Y,Y') \overset{ \text{def} }{=} \hocolim{k>0}
\Omega^{k} k^{\,\Gamma}_{p,k} (Y,Y').
\]
\end{DefRef}

The Fibration Theorem \refT{CPF} can be applied to prove straightforward generalizations of bounded excision.

\begin{DefRef}{THISONE}
The quotient map of categories induces the equivariant map
$K^{\Gamma}_p (Y) \rightarrow
K^{\Gamma}_p (Y,Y')$
and
the map of the fixed points
$K^{\Gamma}_p (Y)^{\Gamma} \rightarrow
K^{\Gamma}_p (Y,Y')^{\Gamma}$.
\end{DefRef}

\begin{PropRef}{FibPair}
In the case of a bounded action of $\Gamma$ on $Y$, there is a homotopy fibration
\[
K^{\Gamma}_p ( Y' )^{\Gamma} \longrightarrow
K^{\Gamma}_p ( Y )^{\Gamma} \longrightarrow K^{\Gamma}_p ( Y,Y' )^{\Gamma}.
\]
\end{PropRef}

\begin{proof}
The inclusion of the subspace $Y'$ induces isomorphisms of categories
$\mathcal{C} (Y') \cong \mathcal{C} (Y)_{<Y'}$ and $\mathcal{C}^{\Gamma}_p ( Y' )^{\Gamma} \cong \mathcal{C}^{\Gamma}_p ( Y )_{<Y'}^{\Gamma}$.
For the second fibration, one should observe that $\mathcal{C}^{\Gamma}_p ( Y )^{\Gamma}$ is
$\mathcal{C}^{\Gamma}_p ( Y )_{<Y'}^{\Gamma}$-filtered.
\end{proof}

\begin{ThmRefName}{BDDEXCI25}{Bounded excision}
If $U_1$ and $U_2$ are a coarsely antithetic pair of subsets of $Y$ which form a cover of $Y$, and the action of $\Gamma$ on $Y$ is bounded, then
\[
\xymatrix{
 K^{\Gamma}_p (U_1 \cap U_2)^{\Gamma} \ar[r] \ar[d]
&K^{\Gamma}_p (U_1)^{\Gamma} \ar[d] \\
 K^{\Gamma}_p (U_2)^{\Gamma} \ar[r]
&K^{\Gamma}_p (Y)^{\Gamma}
}
\]
is a homotopy pushout.
\end{ThmRefName}

\begin{proof}
In view of the isomorphism $\mathcal{C}^{\Gamma}_p (U_1, U_1 \cap U_2) \cong \mathcal{C}^{\Gamma}_p (Y, U_2)$,
we have the weak equivalence
$K^{\Gamma}_p (U_1, U_1 \cap U_2)^{\Gamma} \ \simeq \ K^{\Gamma}_p (Y, U_2)^{\Gamma}$.
We have a map of homotopy fibrations
\[
\begin{CD}
K^{\Gamma}_p (U_1 \cap U_2)^{\Gamma} @>>> K^{\Gamma}_p (U_1)^{\Gamma} @>>> K^{\Gamma}_p (U_1, U_1 \cap U_2)^{\Gamma} \\
@VVV @VVV @VV{\simeq}V \\
K^{\Gamma}_p (U_2)^{\Gamma} @>>> K^{\Gamma}_p (Y)^{\Gamma} @>>> K^{\Gamma}_p (Y, U_2)^{\Gamma}
\end{CD}
\]
which gives the homotopy pushout.
\end{proof}

There are associated relative versions of the excision theorems.

\begin{DefRef}{CoarseInvFun}
Generalizing Definition \refD{THISONE}, if $Y''$ is another coarsely invariant subset of $Y$, then the intersection $Y'' \cap Y'$ is coarsely invariant in both $Y$ and $Y'$, there is an equivariant map
$K^{\Gamma}_p (Y'', Y'' \cap Y')
\rightarrow
K^{\Gamma}_p (Y, Y')$
and on the fixed points
$K^{\Gamma}_p (Y'', Y'' \cap Y')^{\Gamma}
\rightarrow
K^{\Gamma}_p (Y, Y')^{\Gamma}$.
\end{DefRef}

\begin{ThmRefName}{BDDEXCI2}{Relative bounded excision}
Suppose $U_1$, $U_2$, and $Y'$ are three pairwise coarsely antithetic subsets of $Y$
such that
$U_1$ and $U_2$ form a cover of $Y$.
Assuming a bounded action of $\Gamma$ on $Y$,
the commutative square
\[
\xymatrix{
 K^{\Gamma}_{p} (U_1 \cap U_2, Y' \cap U_1 \cap U_2)^{\Gamma} \ar[r] \ar[d]
&K^{\Gamma}_{p} (U_1, Y' \cap U_1)^{\Gamma} \ar[d] \\
 K^{\Gamma}_{p} (U_2, Y' \cap U_2)^{\Gamma} \ar[r]
&K^{\Gamma}_{p} (Y,Y')^{\Gamma} }
\]
induced by inclusions of pairs is
a homotopy pushout.
\end{ThmRefName}

\begin{proof}
This follows from the fact that whenever $C$ is a subset of $Y$ which is coarsely antithetic to $Y'$, the category
$\mathcal{C}^{\Gamma}_p ( Y,Y' )^{\Gamma}$ is
$\mathcal{C}^{\Gamma}_p ( Y,Y' )_{<C}^{\Gamma}$-filtered and $\mathcal{C}^{\Gamma}_p ( Y,Y' )_{<C}^{\Gamma}$ is isomorphic to
$\mathcal{C}^{\Gamma}_p ( C,Y' \cap C )^{\Gamma}$.
The details are left to the reader.
\end{proof}

Then one as easily gets the following corollary.

\begin{ThmRef}{ExcExc}
Suppose that a proper subset $U$ of $Y$ is coarsely invariant under the action of $\Gamma$ and that $U$, $Y - U$, and $Y'$ form a coarse covering by pairwise antithetic subsets.
If the action of $\Gamma$ is bounded, then there is a weak equivalence
\[
K_p^{\Gamma} (Y,Y')^{\Gamma} \simeq K_p^{\Gamma} (Y-U,Y'-U)^{\Gamma}.
\]
\end{ThmRef}

Finally let us formally state the most general version of the excision theorem.
Suppose there is a left action of $\Gamma$ on $Y$ and $\mathcal{U}$ is a coarse covering of $Y$ closed under coarse intersections.

\begin{DefRef}{POIUKI}
The action is $\mathcal{U}$-\textit{bounded} if all coarse families $\{ \mathcal{A} \}$  in $\mathcal{U}$ are closed under the action.  In this case one has the $K$-theory spectra
$K_p^{\Gamma} (\mathcal{A}) = \mathcal{K}_p^{\Gamma} (Y; \mathcal{A})$ and the fixed point spectra $K_p^{\Gamma} (\mathcal{A})^{\Gamma} = \mathcal{K}_p^{\Gamma} (Y; \mathcal{A})^{\Gamma}$.
\end{DefRef}

Of course, an action is $\mathcal{U}$-bounded for any coarse covering $\mathcal{U}$ if the action is by bounded coarse equivalences. 

\begin{ThmRef}{PPPOI2}
Given a $\mathcal{U}$-bounded action on $Y$ for a finite coarse covering $\mathcal{U}$, then
\[
\mathcal{K}_p^{\Gamma} (Y,Y'; \mathcal{U})^{\Gamma} \simeq \hocolim{\mathcal{A} \in \mathcal{U}} K_p^{\Gamma} (\mathcal{A})^{\Gamma} \simeq K_p^{\Gamma} (Y,Y')^{\Gamma}.
\]
\end{ThmRef}

\begin{proof}
Apply the Theorem \refT{BDDEXCI2} inductively to the sets in $\mathcal{U}$.
\end{proof}

It is worth pointing out the connection between the ``core'' of the excision computations here and the main interest of the computations proposed by the Novikov and Borel conjectures, the $K$-theory of the group ring $R[\Gamma]$.

\begin{PropRef}{TPRtriv}
Let $y_0$ be a point in $Y$, then
\[
K^{\Gamma}_{p} (Y)_{<y_0} \simeq K^{\Gamma}_{p} (y_0) = K^{\Gamma,0} (\Gamma)^{\Gamma} \simeq \Knc (R[\Gamma]).
\]
\end{PropRef}

\begin{proof}
	The first equivalence is induced from an isomorphism of bounded categories.  The last equivalence is proved in section VI of \cite{gC:95}.
\end{proof}

\SSecRef{Constructions of bounded actions}{BDDBDD}

Actions by bounded coarse equivalences may seem exceptional, so there is the issue of generating useful examples.

\begin{DefRef}{OrbitMet}
Let $Z$ be any metric space with a free left $\Gamma$-action by isometries.
We assume that the action is properly discontinuous, that is, that for fixed points $z$ and $z'$,
the infimum over $\gamma \in \Gamma$ of the distances $d(z, \gamma z')$ is attained.
Then we define the orbit space metric on $\Gamma \backslash Z$ by
\[
d_{\Gamma \backslash Z} ([z], [z']) = \inf_{\gamma \in \Gamma} d(z, \gamma z').
\]
\end{DefRef}

\begin{LemRef}{JUQASW}
$d_{\Gamma \backslash Z}$ is a metric on $\Gamma \backslash Z$.
\end{LemRef}

\begin{proof}
It is well-known that $d_{\Gamma \backslash Z}$ is a pseudometric.  The fact that $\Gamma$ acts by isometries makes it a metric.
The triangle inequality follows directly from the triangle inequality for $d$.
Symmetry follows from $d (z, \gamma z') = d(\gamma^{-1} z, z') = d(z', \gamma^{-1} z)$.
Finally, $d_{\Gamma \backslash Z} ([z], [z']) = 0$ gives $d (z, \gamma z') = 0$ for some $\gamma \in \Gamma$, so $d (\gamma' z, \gamma' \gamma z') = 0$ for all $\gamma' \in \Gamma$, and so $[z] = [z']$.
\end{proof}

Now suppose $X$ is some metric space with left $\Gamma$-action by isometries.

\begin{DefRef}{Xbdd}
Define
\[
X^{bdd} = X \times_{\Gamma} \Gamma
\]
where
the right-hand copy of $\Gamma$ denotes $\Gamma$ regarded as a metric space with the word-length metric associated to a finite generating set,
the group $\Gamma$ acts by isometries on the metric space $\Gamma$ via left multiplication,
and $X \times_{\Gamma} \Gamma$ denotes the orbit metric space associated to the diagonal left $\Gamma$-action
on $X \times \Gamma$.
We will denote the orbit metric by $d^{bdd}$.
\end{DefRef}

The natural left action of $\Gamma$ on $X^{bdd}$ is given by $\gamma [x,e] = [\gamma x, e]$.

\begin{DefRef}{leftbdd}
A left action of $\Gamma$ on a metric space $X$ is \textit{bounded} if for each element $\gamma \in \Gamma$ there is a number $B_{\gamma} \ge 0$ such that $d(x, \gamma x) \le B_{\gamma}$ for all $x \in X$.
\end{DefRef}

\begin{LemRef}{JUQASW2}
If the left action of $\Gamma$ on a metric space $X$ is bounded, and $B \colon \Gamma \to [0, \infty)$ is a function as above, then there is a real function $B_{\ast} \colon [0, \infty) \to [0, \infty)$ such that $\vert \gamma \vert \le s$ implies $B_{\gamma} \le B_{\ast} (s)$.
\end{LemRef}

\begin{proof}
One simply takes $B_{\ast} (s) = \max \{ B_{\gamma} \mid \vert \gamma \vert \le s \}$.
\end{proof}

\begin{PropRef}{ECEjkiu}
The natural action of $\Gamma$ on $X^{bdd}$ is bounded.
\end{PropRef}

\begin{proof}
If $\vert \gamma \vert = d_{\Gamma} (e, \gamma)$ is the norm in $\Gamma$, we choose $B_{\gamma} = \vert \gamma \vert$.
Now
\begin{align}
d^{bdd} ([x,e], [\gamma x,e]) &= \inf_{\gamma' \in \Gamma} d^{\times} ((x,e), \gamma' (\gamma x,e)) \notag \\
&\le d^{\times} ((x,e), \gamma^{-1} (\gamma x,e)) \notag \\
&= d^{\times} ((x,e), (x,\gamma^{-1})) = d_{\Gamma} (e,\gamma^{-1}) = \vert \gamma^{-1} \vert = \vert \gamma \vert, \notag
\end{align}
where $d^{\times}$ stands for the max metric on the product $X \times \Gamma$.
\end{proof}

\begin{DefRef}{leftbdd44}
Let $b \colon X \to X^{bdd}$ be the natural map given by $b(x) = [x,e]$ in the orbit space $X \times_{\Gamma} \Gamma$.
\end{DefRef}

\begin{PropRef}{HJDSEO}
The map $b \colon X \to X^{bdd}$ is a coarse map.
\end{PropRef}

\begin{proof}
Suppose $d^{bdd} ([x_1,e], [x_2,e]) \le D$, then $d^{\times} ((x_1,e), (\gamma x_2,\gamma)) \le D$ for some $\gamma \in \Gamma$,
so $d(x_1, \gamma x_2) \le D$ and $\vert \gamma \vert \le D$.
Since the left action of $\Gamma$ on $X^{bdd}$ is bounded, there is a function $B_{\ast}$ guaranteed by Lemma \refL{JUQASW2}.
Now
\[
d (x_1, x_2) \le d(x_1, \gamma x_2) + d (x_2, \gamma x_2)
\le D + B_{\ast} (D).
\]
This verifies that $b$ is proper.  It is clearly distance reducing, therefore uniformly expansive with $l(r)=r$.
\end{proof}

If we think of $X^{bdd}$ as the set $X$ with the metric induced from the bijection $b$, the map $b$ becomes the coarse identity map between the metric space $X$ with a left action of $\Gamma$ and the metric space $X^{bdd}$ where the \textit{action is made left-bounded}.

Of course we now have a construction we can apply in any instance of a free properly discontinuous action by $\Gamma$ on a metric space $X$.  In some cases $X$ has a very canonical set of metrics, for example in the case of a fundamental group acting cocompactly on the universal cover $X$.  On the other hand, if $\Gamma$ is the fundamental group of a manifold embedded in a Euclidean space, the normal bundle doesn't have such a set of metrics.  The authors are particularly interested in using left-bounded metrics on the normal bundles.  In these cases the map induced from $b$ preserves the most relevant $K$-theoretic information.

\SecRef{Question, Example, Discussion}{MExD}

\SSecRef{Question}{Question}

We have seen that the lax limit of an action plays an important role in constructing the fixed point spectrum.  If there are two or more natural actions and their comparison is required, one is forced to consider what we will call a \textit{faux lax limit}.  We will then argue that the Karoubi filtration techniques are no longer sufficient to analyze this category.

We exploit the idea that 
the category $\mathcal{C}_p^{\Gamma} (Y,Y')_{\alpha}^{\Gamma}$ and its exact structure can be described independently from the equivariant category $\mathcal{C}_p^{\Gamma} (Y,Y')_{\alpha}$. 

\begin{DefRef}{WGpre2}
First notice that while we do not specify an action of $\Gamma$ on $Y$, there is an action of $\Gamma$ on $\mathcal{C}_{\Gamma} (Y,Y')$ induced from the left multiplication action of the group on itself. 

The \textit{fibred faux lax limit} is a category $\W^{\Gamma} (Y,Y')$ with objects which are sets of data $( \{ F_{\gamma} \},\{ \psi_{\gamma} \} )$ where
\begin{itemize}
\item $F_{\gamma}$ is an object of $\mathcal{C}_{\Gamma} (Y,Y')$ for each $\gamma$ in $\Gamma$,
\item $\psi_{\gamma}$ is an isomorphism $F_e \to F_{\gamma}$ in $\mathcal{C}_{\Gamma} (Y,Y')$,
\item $\psi_{\gamma}$ has filtration $0$ when viewed as a morphism in $\mathcal{C} (\Gamma, \mathcal{C}(Y))$,
\item $\psi_e = \id$,
\item $\psi_{\gamma_1 \gamma_2} = \gamma_1 \psi_{\gamma_2} \circ \psi_{\gamma_1}$ for all $\gamma_1$, $\gamma_2$ in $\Gamma$.
\end{itemize}
The morphisms $( \{ F_{\gamma} \},\{ \psi_{\gamma} \} ) \to ( \{ F'_{\gamma} \},\{ \psi'_{\gamma} \} )$
are collections $ \{ \phi_{\gamma} \} $, where each $\phi_{\gamma}$ is a morphism $F_{\gamma} \to F'_{\gamma}$ in $\mathcal{C}_{\Gamma} (Y,Y')$, such that the squares
\[
\begin{CD}
F_e @>{\psi_{\gamma}}>>  F_{\gamma} \\
@V{\phi_e}VV @VV{ \phi_{\gamma}}V \\
F'_e @>{\psi'_{\gamma}}>> F'_{\gamma}
\end{CD}
\]
commute for all $\gamma \in \Gamma$.
\end{DefRef}

The additive structure on $\W^{\Gamma} (Y,Y')$ is induced from that on $\mathcal{C}_{\Gamma} (Y,Y')$.
For any action $\alpha$ on $Y$ by bounded coarse equivalences, a subset $Y'$ is coarsely invariant, so there is the induced action on the pair $(Y, Y')$.  Now the lax limit
$\mathcal{C}_p^{\Gamma} (Y,Y')^{\Gamma}_{\alpha}$ is an additive subcategory of $\W^{\Gamma} (Y,Y')$.
The embedding $E_{\alpha}$ is realized by sending the object $(F, \psi )$ of $\mathcal{C}_{\Gamma} (Y,Y')^{\Gamma}_{\alpha}$ to
$( \{ \alpha_{\gamma} F \}, \{ \psi (\gamma) \} )$.
On the morphisms, $E_{\alpha} (\phi) = \{ {\alpha}_{\gamma} \phi \}$.

\begin{DefRef}{JKvvv}
The spectrum $w^{\Gamma} (Y,Y')$ is defined as the $K$-theory spectrum of $\W^{\Gamma} (Y,Y')$,
so there are induced map of spectra
$\varepsilon_{\alpha}
\colon
k_p^{\Gamma} (Y,Y')^{\Gamma}_{\alpha}
\rightarrow
w^{\Gamma} (Y,Y')$.
Similarly, there are categories $\W^{\Gamma,k} (Y,Y')$ and the evident exact inclusions
\[
E^k_{\alpha} \colon \mathcal{C}_p^{\Gamma,k} (Y,Y')^{\Gamma}_{\alpha} \longrightarrow \W^{\Gamma,k} (Y,Y').
\]
If the $K$-theory of $\W^{\Gamma,k} (Y,Y')$ is denoted by $w^{\Gamma,k} (Y,Y')$ then
the nonconnective delooping of $w^{\Gamma} (Y,Y')$ can be constructed as
\[
W^{\Gamma} (Y,Y') \, = \, \hocolim{k>0}
\Omega^{k} w^{\Gamma,k} (Y,Y').
\]
Finally, we have induced a map
\[
E_{\alpha} \colon K_p^{\Gamma} (Y,Y')^{\Gamma}_{\alpha} \longrightarrow W^{\Gamma} (Y,Y')
\]
for every choice of the action $\alpha$ on $Y$ by bounded coarse equivalences.
\end{DefRef}

We will argue in the next section that $\W^{\Gamma} (Y,Y')$ can't possess Karoubi filtrations required for a familiar proof of the excision property.  There are other reasons known to us, to be explained in a follow-up paper, why we don't believe the straight-on generalization of the Bounded Excision Theorem \refT{BDDEXCI25} to be true.  

We ask the following question.

\begin{QueRef}{HBVZAQ}
	Under what conditions on the group $\Gamma$ and the ring $R$ does the spectrum valued functor $W^{\Gamma} (Y,Y')$ have the bounded excision property?  Explicitly, under what conditions is the map 
	\[
\rho \colon \mathcal{W}^{\Gamma} (Y; \mathcal{U})^{\Gamma} = \hocolim{\mathcal{A} \in \mathcal{U}} W^{\Gamma} (Y,Y')^{\Gamma}_{\mathcal{A}} \longrightarrow W^{\Gamma} (Y,Y')^{\Gamma}
\]
a weak equivalence for a satuarated antithetic covering $\mathcal{U}$?
\end{QueRef}

Our experience from \cite{gCbG:15} leads us to the following guess.

\begin{ConjRef}{MORCAQ}
	The map $\rho$ in Question \refQ{HBVZAQ} is a weak equivalence whenever $\Gamma$ has a coarse embedding in Hilbert space and $R$ is a regular Noetherian ring with finite global dimension.
\end{ConjRef}

\SSecRef{Example}{Example}

In order to illustrate at this juncture the difficulties with the Karoubi filtrations that arise, we will construct and analyze in this section a simplified and absolute version $W (Y)$ which deals with just the fibre but has the essential features of $W^{\Gamma} (Y,Y')$.
Just as before, the category $\W (Y)$ is not 
the lax limit of any specific $\Gamma$-action.  

\begin{DefRef}{WGpre2or}
The category $\mathcal{W} (Y)$ is what we call a \textit{faux lax limit}.  In this example, it has objects which are sets of data $F = ( \{ F_{\gamma} \},\{ \psi_{\gamma} \} )$ where
\begin{itemize}
\item $F_{\gamma}$ is an object of $\mathcal{C} (Y,R)$ for each $\gamma$ in $\Gamma$,
\item $\psi_{\gamma}$ is an isomorphism $F_e \to F_{\gamma}$ in $\mathcal{C} (Y,R)$,
\item $\psi_e = \id$.
\end{itemize}
A morphism $f \colon F \to F'$
is a collection of morphisms $ \{ f_{\gamma} \colon F_{\gamma} \to F'_{\gamma} \}$ in $\mathcal{C} (Y,R)$, such that the squares
\[
\begin{CD}
F_e @>{\psi_{\gamma}}>>  F_{\gamma} \\
@V{f_e}VV @VV{ f_{\gamma}}V \\
F'_e @>{\psi'_{\gamma}}>> F'_{\gamma}
\end{CD}
\]
commute for all $\gamma \in \Gamma$.
\end{DefRef}

The additive structure on $\mathcal{W} (Y)$ is induced from that on $\mathcal{C} (Y,R)$. 
What this means is that the operation $\oplus$ has the property $(F \oplus G)_{\gamma} (m) =F_{\gamma} (m) \oplus G_{\gamma} (m)$. 
Notice also that in order to define a morphism $f$ it suffices to give one component $f_0 \colon F_0 \to F'_0$.  All other components $f_{\gamma}$ are enforced by the virtue of $\psi_{\gamma}$ and $\psi'_{\gamma}$ being isomorphisms.

The metric space $Y$ is the integers $\mathbb{Z}$ with the standard metric invariant under translations.
The group in this example is the additive group of the integers, the infinite cyclic group, which is denoted by $C$.
The essence of failure of the Karoubi filtration technique can be seen already here.

We want to demonstrate an object of the faux lax limit $\mathcal{W} (\mathbb{Z})$ which cannot have $\oplus$-decompositions that are part of a Karoubi filtration.  For the filtering subcategory we choose $\mathcal{W} (\mathbb{Z})_{< \mathbb{Z}^{\le 0}}$.
The object $G$ is defined by $G_n = \bigoplus_{\mathbb{Z}} R$.  The isomorphism $\psi_n \colon G_0 \to G_n$ is given by the bi-infinite matrix over $R$ in the standard basis, with rows and columns enumerated from left to right and from top down respectively.  To construct a matrix for $\psi_n$ one starts with the infinite identity matrix and replaces the square block $[0,n] \times [-n,0]$ with the lower triangular matrix of $1$s below the diagonal.  In other words, $\psi_n$ restricts to the identity on both submodules $G_0 ((- \infty, -1])$ and $G_0 ([n+1, + \infty))$.  It is given by the lower triangular matrix of $1$s on $G_0 ([0,n])$, which is a nonsingular matrix.
The whole $\psi_n$ is an isomorphism bounded by $n+1$. 
The important feature of $\psi_n$ is that the submodule $(G_{0})_0 = R$ maps onto the direct summand in $G_n ([0,n])$ generated by the element $(1, 1, \ldots, 1)$.  This direct summand, however, is not a direct summand in the Pedersen-Weibel object $G_n$ in $\mathcal{C} (\mathbb{Z})$.  If we now restrict $\psi_n$ to $G_0 ((- \infty, 0])$ then the minimal direct summand of $G_n$ containing the image will be the submodule $G_n  ((- \infty, n])$.  We remind the reader that the decompositions in the definition of Karoubi filtrations $E_{\alpha} \oplus D_{\alpha}$ in our case are the direct sum decompositions in $\mathcal{W} (\mathbb{Z})$.  On the level of each element $n$ of $C$, for all $k$ in $\mathbb{Z}$ this enforces the direct sum $G_n (k) = (E_{\alpha})_n (k) \oplus (D_{\alpha})_n (k)$.

In order to see that no proper decompositions can give a Karoubi filtration of $\mathcal{W} (\mathbb{Z})$ by $\mathcal{W} (\mathbb{Z})_{< \mathbb{Z}^{\le 0}}$, we define an object
$F$ by $F_n = G_0 ((- \infty, 0])$ and assigning all structural isomorphisms to be the identity maps.  It is clearly an object of $\mathcal{W} (\mathbb{Z})_{< \mathbb{Z}^{\le 0}}$.  
The morphism $f \colon F \to G$ is given by making $f_n$ the restriction $\psi_n$ to the submodule $G_0 ((- \infty, 0])$. 
On the level of $f_0 \colon F_0 \to G_0$, this morphism is specified by the inclusion $G_0 ((- \infty, 0]) \subset G_0$.  We claim that $f$ fails to factor through any proper decomposition of $G$.

Indeed, we saw that for each $n$ there is a minimal direct summand $G_n  ((- \infty, n])$ that should be a submodule of such $(E_{\alpha})_n$.  There is also the stalk $G_0 (n)$ that maps identically onto $G_n (n)$.  Since $E_{\alpha}$ has to have the structural isomorphisms which are restrictions of $\psi_n$, we conclude that for all $n$ the whole submodule $G_0 (n)$ must be contained in $(E_{\alpha})_0$.  This enforces $(E_{\alpha})_0 = G_0$ and therefore $E_{\alpha} = G$.

\SSecRef{Discussion}{Discussion}
One may wonder about the options available to repair the rigidity of the $\oplus$-structure in $\mathcal{W} (\mathbb{Z})$ in order to accommodate images of maps such as $f$.  This would require one to extend the $\oplus$ operation to direct sums of modules which are not themselves Pedersen-Weibel objects even though their sums might be.  To do this coherently, one has to include the filtered free $R$-modules $F$ which assign to each subset $S$ of ${\mathbb{Z}}$, or more generally a proper metric space $M$, a free $R$-submodule $F(S)$ of $F$.  It may be necessary to require that if $S$ is a subset of $T$ then $F(S)$ is a direct summand of $F(T)$.  But unlike the situation in the Pedersen-Weibel objects, one should not expect a splitting $F(T) = F(S) \oplus F(T \setminus S)$.  This means that one should not expect the $K$-theory of this new additive category to be the bounded $K$-theory of $M$. A viable theory of this kind can be constructed and, in fact, does have a successful analogue in $L$-theory.  
                                
Another resolution with good excision properties can be built as the fibred bounded $G$-theory.  This is a theory generalizing \cite{gCbG:00} where all desired localization and excision properties (true and false) of bounded $K$-theory do hold without any conditions on the group and the ring, as long as the coefficient ring $R$ is Noetherian.  In many ways this development mimics the classical relationship between $K$-theory and $G$-theory of Noetherian rings.


\begin{thebibliography}{99}

\bibitem{sA:16}
S. Arnt, \textit{Fibred coarse embeddability of box spaces and proper isometric affine actions on $L^p$ spaces}, Bull. Belg. Math. Soc. Simon Stevin \textbf{23} (2016), 21--32.

\bibitem{sBbG:18}
{S. Beckhardt and B.~Goldfarb}, 
\textit{Extension properties of groups with asymptotic property C and straight finite decomposition complexity}, Topol. Appl. \textbf{239} (2018), 181--190.

\bibitem{mCeP:97}
{M. Cardenas and E.K. Pedersen},
{\it On the Karoubi filtration of a category},
$K$-theory, {\bf 12} (1997), 165--191.

\bibitem{gC:95}
{G. Carlsson},
{\it Bounded $K$-theory and the assembly map
in algebraic $K$-theory},
in {\it Novikov conjectures, index theory and rigidity},
{\it Vol. 2}
(S.C. Ferry, A. Ranicki, and J. Rosenberg, eds.),
Cambridge U. Press (1995), 5--127.

\bibitem{gCbG:03}
{G.~Carlsson and B.~Goldfarb}, \textit{On homological coherence of
discrete groups}, J.~Algebra {\bf 276} (2004), 502--514.

\bibitem{gCbG:04}
\bysame, \textit{The integral K-theoretic Novikov conjecture for
groups with finite asymptotic dimension}, Invent. Math. {\bf
157} (2004), 405--418.

\bibitem{gCbG:00}
\bysame,
\textit{Controlled algebraic $G$-theory, I}, J. Homotopy Relat. Struct. \textbf{6} (2011), 119--159.

\bibitem{gCbG:13}
\bysame, \textit{Algebraic $K$-theory of geometric groups}, preprint, 2015.
\texttt{arXiv:1305.3349}

\bibitem{gCbG:14}
\bysame, \textit{Equivariant stable homotopy methods in the algebraic K-theory of infinite groups}, preprint, 2014.
\texttt{arXiv:1412.3483}

\bibitem{gCbG:15}
\bysame, \textit{On modules over infinite group rings}, Int. J. Algebra Comput. \textbf{26} (2016), 1--16.


\bibitem{gCeP:93}
{G. Carlsson and E.K. Pedersen}, {\it Controlled algebra and the
Novikov conjecture for K- and L-theory}, Topology {\bf 34}
(1993), 731--758.

\bibitem{xCqWxW:13}
X. Chen, Q. Wang, and X. Wang, \textit{Characterization of the Haagerup property by fibred coarse embedding into Hilbert space}, Bull. London Math. Soc. \textbf{45} (2013), 1091--1099.

\bibitem{xCqWgY:13}
X. Chen, Q. Wang, and G. Yu,
\textit{The maximal coarse Baum-Connes conjecture for spaces which admit a fibred coarse embedding into Hilbert space}, Adv. Math. \textbf{249} (2013), 88--130.

\bibitem{mF:14}
M. Finn-Sell, 
\textit{Fibred coarse embeddings, a-T-menability and the coarse analogue of the Novikov conjecture}, 
J. Funct. Anal. \textbf{267} (2014), 3758--3782.

\bibitem{bGjG:18}
{B. Goldfarb and J. Grossman}, 
\textit{Coarse coherence of metric spaces and groups and its permanence properties}, preprint, 2018. \texttt{arXiv:1804.00944}

\bibitem{dK:15}
{D. Kasprowski}, 
\textit{On the K-theory of groups with finite decomposition complexity},
Proc. London Math. Soc. \textbf{120} (2015), 565--592.

\bibitem{mMhS:18}
M. Mimura and H. Sako, 
\textit{Group approximation in Cayley topology and coarse geometry, Part II: Fibered coarse embeddings},
preprint, 2018.
\texttt{arXiv:1804.10614}

\bibitem{eP:84}
{E.K. Pedersen},
{\it On the $K_{-i}$-functors},
J. Algebra \textbf{90} (1984), 461--475.

\bibitem{ePcW:85}
{E.K. Pedersen and C. Weibel},
{\it A nonconnective delooping of algebraic $K$-theory},
in {\it Algebraic and geometric topology}
(A. Ranicki, N. Levitt, and F. Quinn, eds.),
Lecture Notes in Mathematics {\bf 1126},
Springer-Verlag (1985), 166--181.

\bibitem{ePcW:89}
{\bysame}, {\it $K$-theory homology of spaces}, in {\it Algebraic
topology} (G.~Carlsson, R.L.~Cohen, H.R.~Miller, and
D.C.~Ravenel, eds.), Lecture Notes in Mathematics {\bf 1370},
Springer-Verlag (1989), 346--361.


\bibitem{jR:03}
{J. Roe},
{\it Lectures on coarse geometry}, University Lecture Series, vol.~31, American Mathematical Society, 2003.

\bibitem{rT:83}
{R.W. Thomason},
{\it The homotopy limit problem},
in {\it Proceedings of the Northwestern homotopy theory conference}
(H.R. Miller and S.B. Priddy, eds.),
Contemp. Math. {\bf 19} (1983), 407--420.


\bibitem{kW:10}
{K.~Whyte}, \textit{Coarse bundles}, preprint, 2010. \texttt{arXiv:1006.3347}

\bibitem{gY:98}
{G.~Yu}, \textit{The Novikov conjecture for groups with finite asymptotic dimension}, Annals of Math. \textbf{147} (1998), 325--355.

\bibitem{gY:00}
\bysame, \textit{The coarse Baum-Connes conjecture for spaces which admit a uniform embedding into Hilbert space}, Invent. Math. \textbf{139} (2000), 201--240 

\end{thebibliography}
\end{document}